# Radial Basis Function Approximations: Comparison and Applications

Zuzana Majdisova[a,*], Vaclav Skala[a]

[a] *Department of Computer Science and Engineering, Faculty of Applied Sciences, University of West Bohemia, Univerzitni 8, CZ 30614 Plzen, Czech Republic*

**Abstract**

Approximation of scattered data is often a task in many engineering problems. The Radial Basis Function (RBF) approximation is appropriate for large scattered (unordered) datasets in $d$-dimensional space. This approach is useful for a higher dimension $d > 2$, because the other methods require the conversion of a scattered dataset to an ordered dataset (i.e. a semi-regular mesh is obtained by using some tessellation techniques), which is computationally expensive. The RBF approximation is non-separable, as it is based on the distance between two points. This method leads to a solution of Linear System of Equations (LSE) $\mathbf{Ac} = \mathbf{h}$.

In this paper several RBF approximation methods are briefly introduced and a comparison of those is made with respect to the stability and accuracy of computation. The proposed RBF approximation offers lower memory requirements and better quality of approximation.

## 1. Introduction

Radial Basis Functions (RBFs) are widely used across of many fields solving technical and non-technical problems. A RBF method was originally introduced by [1] and it is an effective tool for solving partial differential equations in engineering and sciences. Moreover, RBF applications can be found in neural networks, fuzzy systems, pattern recognition, data visualization, medical applications, surface reconstruction [2], [3], [4], [5], reconstruction of corrupted images [6], [7], etc. The RBF approximation technique is really meshless and is based on collocation in a set of scattered nodes. This method is independent with respect to the dimension of the space. The computational cost of RBF approximation increases nonlinearly with the number of points in the given dataset and linearly with the dimensionality of data.

*Corresponding author
*Email address:* majdisz@kiv.zcu.cz (Zuzana Majdisova)
*URL:* www.vaclavskala.eu (Vaclav Skala)

There are two main groups of basis functions: global RBFs and Compactly Supported RBFs (CS−RBFs) [8]. Fitting scattered data with CS−RBFs leads to a simpler and faster computation, because a system of linear equations has a sparse matrix. However, approximation using CS−RBFs is quite sensitive to the density of approximated scattered data and to the choice of a shape parameter. Global RBFs lead to a linear system of equations with a dense matrix and their usage is based on sophisticated techniques such as the fast multipole method [9]. Global RBFs are useful in repairing incomplete datasets and they are insensitive to the density of approximated data.

## 2. RBF Approximation using Lagrange Multipliers

RBF approximation introduced by Fasshauer [10] (Chapter 19) is based on Lagrange multipliers. In this section, the properties of this method will be briefly summarized.

This RBF approximation is formulated as a constrained quadratic optimization problem. The goal of this method is to approximate the given dataset by function:

$$f(\mathbf{x}) = \sum_{j=1}^{M} c_j \phi(\|\mathbf{x} - \boldsymbol{\xi}_j\|), \tag{1}$$

where the approximating function $f(\mathbf{x})$ is represented as a sum of $M$ RBFs, each associated with a different reference point $\boldsymbol{\xi}_j$, and weighted by an appropriate coefficient $c_j$. Therefore, it is necessary to determine the vector of weights $\mathbf{c} = (c_1, \ldots, c_M)^T$, which leads to the minimization of the quadratic form:

$$\frac{1}{2}\mathbf{c}^T \mathbf{Q} \mathbf{c}, \tag{2}$$

where $\mathbf{Q}$ is some $M \times M$ symmetric positive definite matrix. This quadratic form is minimized subject to the $N$ linear constraints $\mathbf{Ac} = \mathbf{h}$, where $\mathbf{A}$ is an $N \times M$ matrix with full rank, and the right-hand side $\mathbf{h} = (h_1, \ldots, h_N)^T$ is given. Thus the constrained quadratic minimization problem can be described as a LSE:

$$F(\mathbf{c}, \boldsymbol{\lambda}) = \frac{1}{2}\mathbf{c}^T\mathbf{Q}\mathbf{c} - \boldsymbol{\lambda}^T(\mathbf{A}\mathbf{c} - \mathbf{h}), \tag{3}$$

where $\boldsymbol{\lambda} = (\lambda_1, \ldots, \lambda_N)^T$ is the vector of Lagrange multipliers, and we need to find the minimum of (3) with respect to $\mathbf{c}$ and $\boldsymbol{\lambda}$. This leads to solving the following system:

$$\begin{aligned} \frac{\partial F(\mathbf{c}, \boldsymbol{\lambda})}{\partial \mathbf{c}} &= \mathbf{Q}\mathbf{c} - \mathbf{A}^T\boldsymbol{\lambda} = \mathbf{0} \\ \frac{\partial F(\mathbf{c}, \boldsymbol{\lambda})}{\partial \boldsymbol{\lambda}} &= \mathbf{A}\mathbf{c} - \mathbf{h} = \mathbf{0} \end{aligned} \tag{4}$$

or, in matrix form:

$$\begin{pmatrix} \mathbf{Q} & -\mathbf{A}^T \\ \mathbf{A} & \mathbf{0} \end{pmatrix} \begin{pmatrix} \mathbf{c} \\ \boldsymbol{\lambda} \end{pmatrix} = \begin{pmatrix} \mathbf{0} \\ \mathbf{h} \end{pmatrix}, \tag{5}$$

where $Q_{i,j} = \phi(\|\boldsymbol{\xi}_i - \boldsymbol{\xi}_j\|)$ and $\mathbf{Q}$ is a symmetric matrix. Equation (5) is then solved.

It should be noted that we want to minimize $M$ in order to reduce the computational cost of the approximated value $f(\mathbf{x})$ as much as possible.

## 3. RBF Approximation

Another approach is RBF interpolation, which is based on a solution of a linear system of equations (LSE) [11]:

$$\mathbf{A}\mathbf{c} = \mathbf{h}, \tag{6}$$

where $\mathbf{A}$ is a matrix of this system, $\mathbf{c}$ is a column vector of variables and $\mathbf{h}$ is a column vector containing the right sides of equations. In this case, $\mathbf{A}$ is an $N \times N$ matrix, where $N$ is the number of given points, the variables are weights for basis functions and the right sides of equations are values in the given points. The disadvantage of RBF interpolation is the large and usually ill-conditioned matrix of the LSE. Moreover, in the case of an oversampled dataset or intended reduction, we want to reduce the given problem, i.e. reduce the number of weights and used basis functions, and preserve good precision of the approximated solution. The approach, which includes the reduction, is called RBF approximation. In the following, the method recently introduced in [11] is described in detail.

For simplicity, we assume that we have an unordered dataset $\{\mathbf{x}_i\}_1^N$ in $E^2$. However, note that this approach is generally applicable for $d$-dimensional space. Further, each point $\mathbf{x}_i$ from the dataset is associated with vector $\mathbf{h}_i \in E^p$ of the given values, where $p$ is the dimension of the vector, or scalar value $h_i \in E^1$. For an explanation of the RBF approximation, let us consider the case when each point $\mathbf{x}_i$ is associated with scalar value $h_i$. Now we extend the given dataset by a set of new reference points $\{\boldsymbol{\xi}_j\}_1^M$, see Fig. 1.

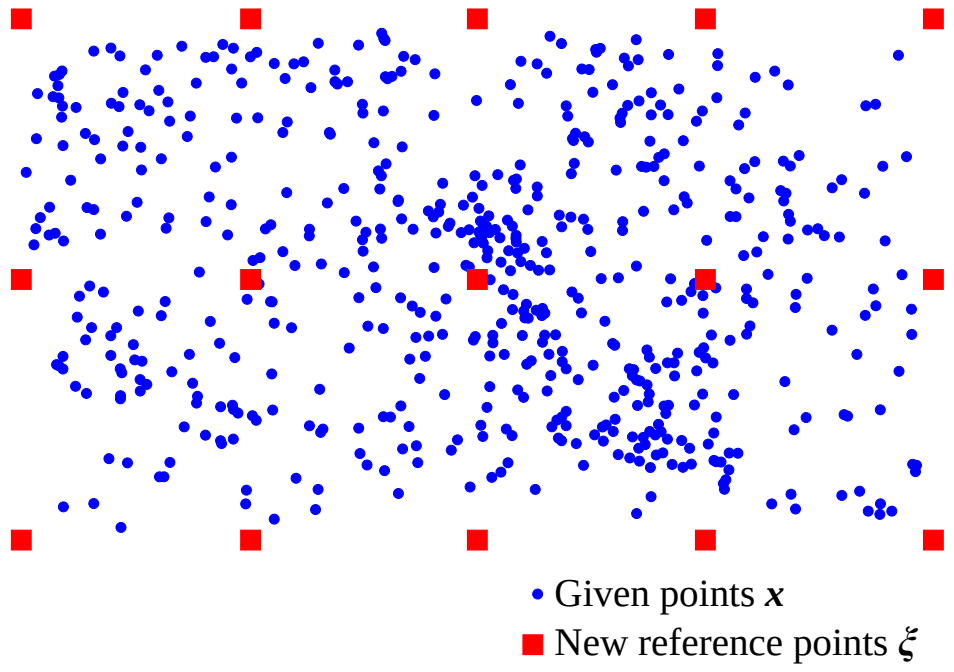

Figure 1: RBF approximation and reduction of points.

These reference points may not necessarily be in a uniform grid. It is appropriate, that their placement reflects the given surface as well as possible. A good placement of the reference points improves the approximation of the underlying data. For example, when a terrain is to be approximated, placement along features such as break lines leads to better approximation results. The number of added reference points $\boldsymbol{\xi}_j$ is $M$, where $M \ll N$. The RBF approximation is based on computing the distance of given point $\mathbf{x}_i$ and reference point $\boldsymbol{\xi}_j$ from the extended dataset.

The approximated value can be determined similarly as for interpolation (see [11]):

$$f(\mathbf{x}) = \sum_{j=1}^{M} c_j \phi(r_j) = \sum_{j=1}^{M} c_j \phi(\|\mathbf{x} - \boldsymbol{\xi}_j\|), \tag{7}$$

where the approximating function $f(\mathbf{x})$ is represented as a sum of $M$ RBFs, each associated with a different reference point $\boldsymbol{\xi}_j$, and weighted by an appropriate coefficient $c_j$.

It can be seen that we get an overdetermined LSE for the given dataset:

$$h_i = f(\mathbf{x}_i) = \sum_{j=1}^{M} c_j \phi(\|\mathbf{x}_i - \boldsymbol{\xi}_j\|) = \sum_{j=1}^{M} c_j \phi_{i,j} \qquad i = 1, \ldots, N. \tag{8}$$

The linear system of equations (8) can be represented as the matrix equation:

$$\mathbf{Ac} = \mathbf{h}, \tag{9}$$

where the number of rows is $N \gg M$ and $M$ is the number of unknown weights $[c_1, \ldots, c_M]^T$, i.e. the number of reference points. Equation (9) can be expressed in the form:

$$\begin{pmatrix} \phi_{1,1} & \cdots & \phi_{1,M} \\ \vdots & \ddots & \vdots \\ \phi_{i,1} & \cdots & \phi_{i,M} \\ \vdots & \ddots & \vdots \\ \phi_{N,1} & \cdots & \phi_{N,M} \end{pmatrix} \begin{pmatrix} c_1 \\ \vdots \\ c_M \end{pmatrix} = \begin{pmatrix} h_1 \\ \vdots \\ h_i \\ \vdots \\ h_N \end{pmatrix}. \tag{10}$$

Thus the presented system is overdetermined, i.e. the number of equations $N$ is higher than number of variables $M$. This LSE can be solved by the least squares method as $\mathbf{A}^T\mathbf{Ac} = \mathbf{A}^T\mathbf{h}$ or singular value decomposition, etc.

## 4. RBF Approximation with Polynomial Reproduction

The method which was introduced in Sect. 3 can theoretically have problems with stability and solvability. Therefore, the RBF approximant (7) is usually extended by polynomial function $P_k(\mathbf{x})$ of degree $k$. Now, the approximated value can be expressed in the form:

$$f(\mathbf{x}) = \sum_{j=1}^{M} c_j \phi(\|\mathbf{x} - \boldsymbol{\xi}_j\|) + P_k(\mathbf{x}). \tag{11}$$

where $\boldsymbol{\xi}_j$ are reference points specified by a user. This leads to solving the LSE:

$$\begin{aligned} h_i = f(\mathbf{x}_i) &= \sum_{j=1}^{M} c_j \phi(\|\mathbf{x}_i - \boldsymbol{\xi}_j\|) + P_k(\mathbf{x}_i) \\ &= \sum_{j=1}^{M} c_j \phi_{i,j} + P_k(\mathbf{x}_i) \qquad i = 1, \ldots, N. \end{aligned} \tag{12}$$

In practice, a linear polynomial:

$$P_1(\mathbf{x}) = \mathbf{a}^T\mathbf{x} + a_0 \tag{13}$$

is used. Geometrically, the coefficient $a_0$ determines the placement of the hyperplane and the expression $\mathbf{a}^T\mathbf{x}$ represents the inclination of the hyperplane.

It can be seen that for $d$-dimensional space a linear system of $N$ equations in $(M+d+1)$ variables has to be solved, where $N$ is the number of points in the given dataset, $M$ is the number of reference points and $d$ is the dimensionality of space, e.g. for $d = 2$ vectors $\mathbf{x}_i$ and $\mathbf{a}$ are given as $\mathbf{x}_i = (x_i, y_i)^T$ and $\mathbf{a} = (a_x, a_y)^T$. Using the matrix notation, we can write for $E^2$:

$$\begin{pmatrix} \phi_{1,1} & \cdots & \phi_{1,M} & x_1 & y_1 & 1 \\ \vdots & \ddots & \vdots & \vdots & \vdots & \vdots \\ \phi_{i,1} & \cdots & \phi_{i,M} & x_i & y_i & 1 \\ \vdots & \ddots & \vdots & \vdots & \vdots & \vdots \\ \phi_{N,1} & \cdots & \phi_{N,M} & x_N & y_N & 1 \end{pmatrix} \begin{pmatrix} c_1 \\ \vdots \\ c_M \\ a_x \\ a_y \\ a_0 \end{pmatrix} = \begin{pmatrix} h_1 \\ \vdots \\ h_i \\ \vdots \\ h_N \end{pmatrix}. \tag{14}$$

Equation (14) can also be expressed in the form:

$$\begin{pmatrix} \mathbf{A} & \mathbf{P} \end{pmatrix} \begin{pmatrix} \mathbf{c} \\ \mathbf{a} \\ a_0 \end{pmatrix} = \mathbf{h}. \tag{15}$$

It can be seen that for $E^2$ we have a linear system of $N$ equations in $(M + 3)$ variables, where $M \ll N$. Thus the presented system is overdetermined again and can also be solved by the method of least squares or singular value decomposition.

## 5. Experimental Results

The above presented methods of the RBF approximation have been tested on synthetic and real datasets. Moreover, different global radial basis functions with shape parameter $\alpha$, see Table 1, and different sets of reference points have been used for testing. These sets of reference points have different types of distributions described in Sect. 5.1.

### *5.1. Distribution of Reference Points*

For these experiments, the following sets of reference points were used:

**Points on regular grid**

This set contains the points on a regular grid in $E^2$.

Table 1: Used global RBFs ($\alpha$ is a shape parameter)

| RBF | $\phi(\mathbf{r})$ |
|---|---|
| Gauss function [12] | $e^{-(\alpha r)^2}$ |
| Inverse Quadric (IQ) | $\frac{1}{1+(\alpha r)^2}$ |
| Thin-Plate Spline (TPS) [13] | $(\alpha r)^2 \log(\alpha r)$ |

**Epsilon points**
This distribution of reference points is described in the following text.

**Epsilon points + AABB corners**
This set of points is determined in the same manner as the previous case. Moreover, the corners of axis aligned bounding box (AABB) of Epsilon points are added to the set of reference points.

**Halton points**
This distribution of points is described in the following text in detail. However, note that this set of reference points equals the subset of the given dataset, for which we determine the RBF approximation.

**Halton points + AABB corners**
This set of reference points is determined in the same manner as Halton points. Moreover, the corners of AABB are added to this set.

#### *5.1.1. Halton points*

Construction of a Halton sequence is based on a deterministic method. This sequence generates well-spaced "draws" points from the interval $[0, 1]$. The sequence uses a prime number as its base and is constructed based on finer and finer prime-based divisions of sub-intervals of the unit interval. The Halton sequence [10] can be described by the following recurrence formula:

$$Halton(p)_k = \sum_{i=0}^{\lfloor \log_p k \rfloor} \frac{1}{p^{i+1}} \left( \left\lfloor \frac{k}{p^i} \right\rfloor \bmod p \right), \tag{16}$$

where $p$ is the prime number and $k$ is the index of the calculated element.

For the $E^2$ space, subsequent prime numbers are used as a base. In this test, $\{2, 3\}$ were used for the Halton sequence and the following sequence of points in a rectangle $(a, b)$ was derived:

$$Halton(2,3) = \Bigg\{ \left(\frac{1}{2}a, \frac{1}{3}b\right), \left(\frac{1}{4}a, \frac{2}{3}b\right), \left(\frac{3}{4}a, \frac{1}{9}b\right), \left(\frac{1}{8}a, \frac{4}{9}b\right), \left(\frac{5}{8}a, \frac{7}{9}b\right), \left(\frac{3}{8}a, \frac{2}{9}b\right), \left(\frac{7}{8}a, \frac{5}{9}b\right), \left(\frac{1}{16}a, \frac{8}{9}b\right), \left(\frac{9}{16}a, \frac{1}{27}b\right), \dots \Bigg\}, \tag{17}$$

where $a$ is the width of the rectangle and $b$ is the height of the rectangle.

Visualization of the dataset with $10^3$ points of the Halton sequence from (17) can be seen in Fig. 2. We can see that the Halton sequence in $E^2$ space covers this space more evenly than randomly distributed uniform points in the same rectangle.

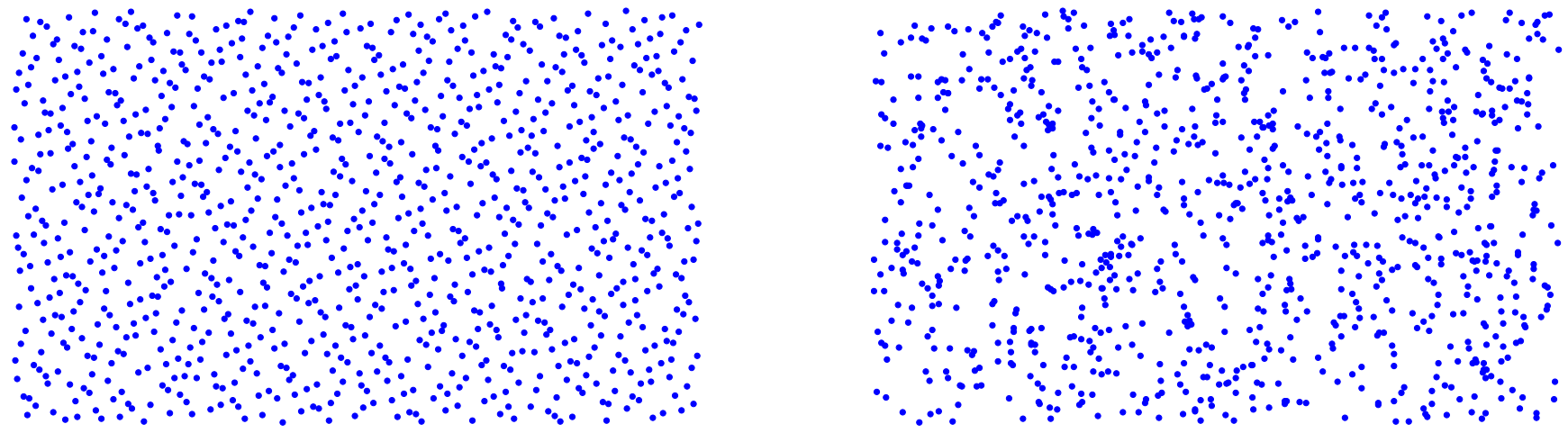

Figure 2: Halton points in $E^2$ generated by $Halton(2, 3)$ (left) and random points in a rectangle with uniform distribution (right). The number of points is $10^3$ in both cases.

### *5.1.2. Epsilon points*

This is a special distribution of points in $E^2$, which is based on a regular grid. Each point is determined as follows:

$$\begin{aligned} P_{ij} &= \Big[ i \cdot \Delta x + \text{rand}(-\varepsilon_x, \varepsilon_x), j \cdot \Delta y + \text{rand}(-\varepsilon_y, \varepsilon_y) \Big], \\ &\varepsilon_x \approx 0.25 \cdot \Delta \text{x}, \qquad i = 0, \dots, N_x, \\ &\varepsilon_y \approx 0.25 \cdot \Delta \text{y}, \qquad j = 0, \dots, N_y, \end{aligned} \tag{18}$$

where $\Delta x$ and $\Delta y$ are real numbers representing the grid spacing, $N_x$ indicates the number of grid columns, $N_y$ is the number of grid rows and $\text{rand}(-\varepsilon_x, \varepsilon_x)$ or rand$(-\varepsilon_y, \varepsilon_y)$ is a random drift with a uniform distribution from $-\varepsilon_x$ to $\varepsilon_x$ or from $-\varepsilon_y$ to $\varepsilon_y$.

Figure 3 presents the dataset with $40 \times 25$, (i.e. $10^3$) epsilon points. Moreover, we can see the comparison of this distribution of points with points on a regular grid.

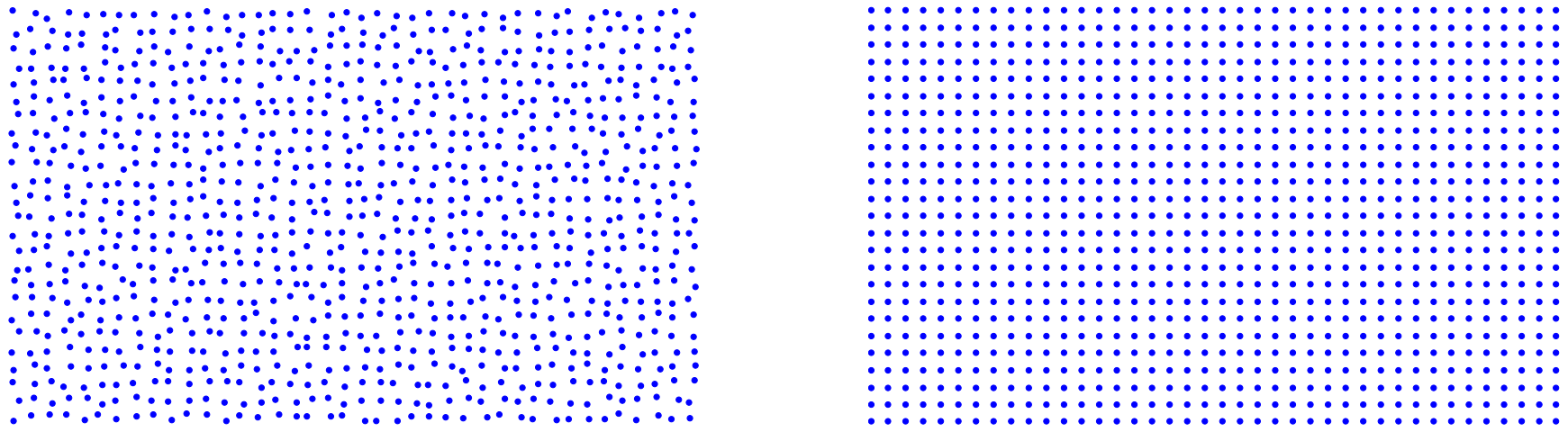

Figure 3: Epsilon points (left) and points on a $2D$ regular grid (right). The number of points is $40 \times 25 = 10^3$ in both cases.

### *5.2. Synthetic Datasets*

The Halton distribution of points was used for synthetic data. Moreover, each point from this dataset is associated with a function value at this point. For this purpose, different functions have been used for experiments. Results for a $2D$ sinc function, see Fig. 4, are presented in this paper.

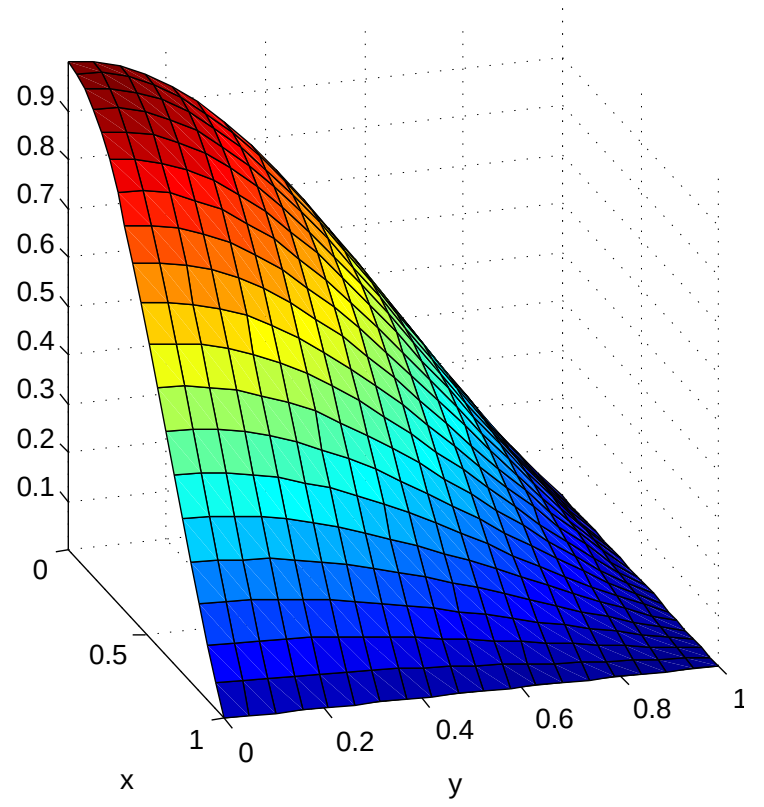

Figure 4: 2D sinc function $\mathrm{sinc}(\pi x) \cdot \mathrm{sinc}(\pi y)$ whose domain is restricted to $[0, 1] \times [0, 1]$.

#### *5.2.1. Examples of RBF Approximation Results*

Some examples of RBF approximation to 1089 Halton data points sampled from a $2D$ sinc function, for a Halton set of reference points, which consists of 81 points, and different RBFs are shown in Fig. 5 and Fig. 6.

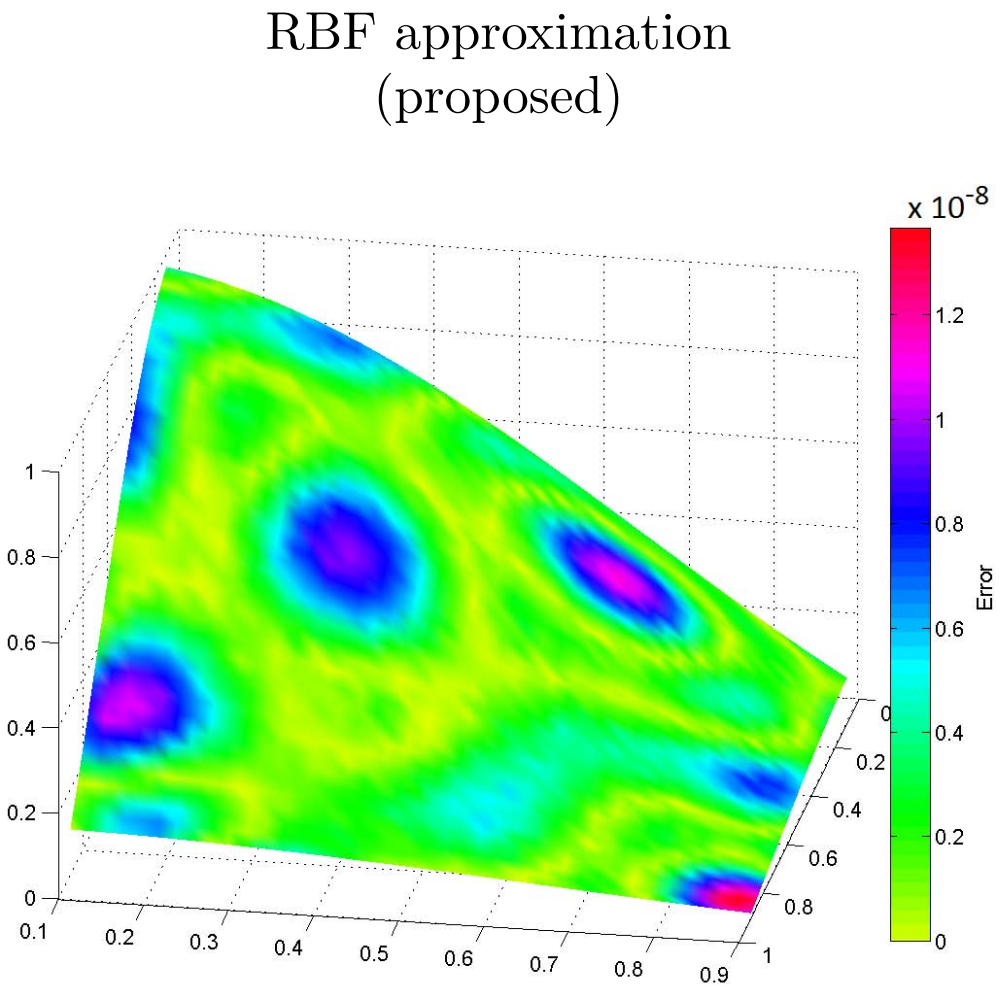

(a) Gauss, $\alpha = 1$, Halton points

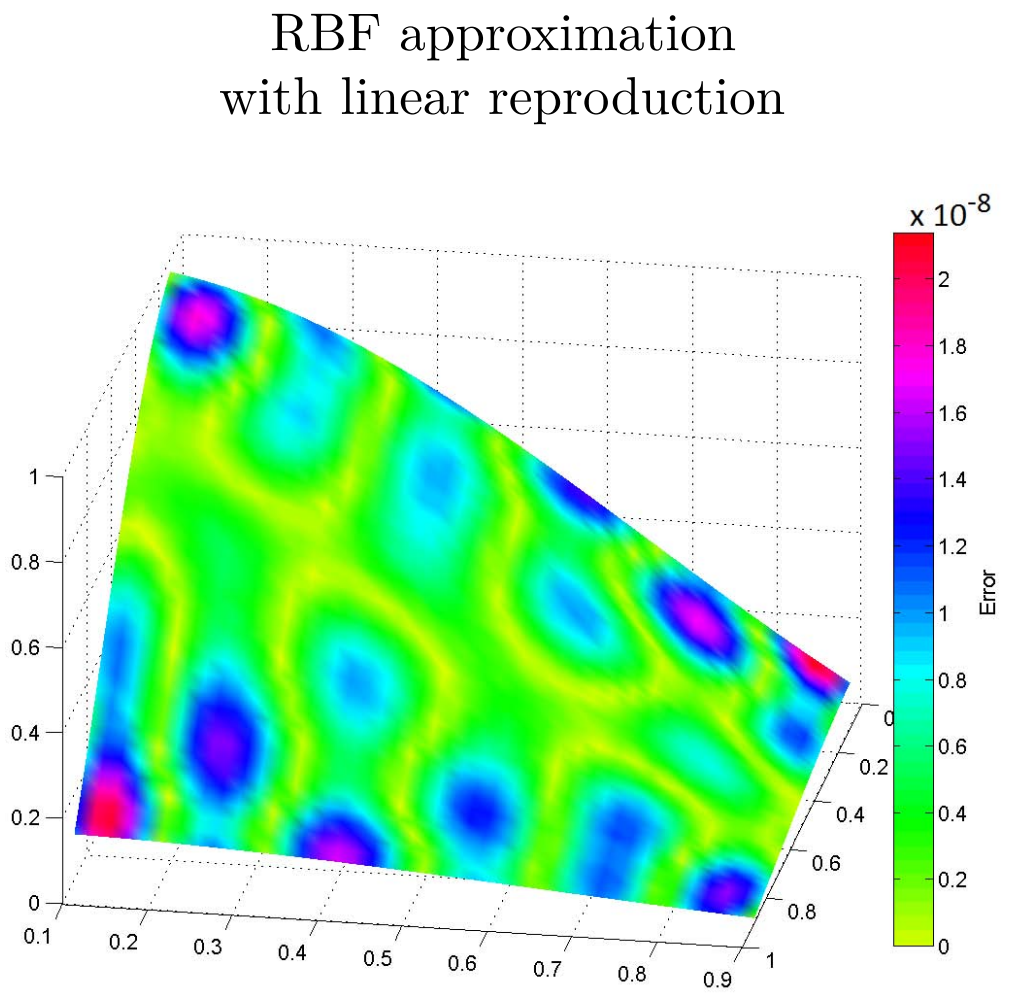

(b) Gauss, $\alpha = 1$, Halton points

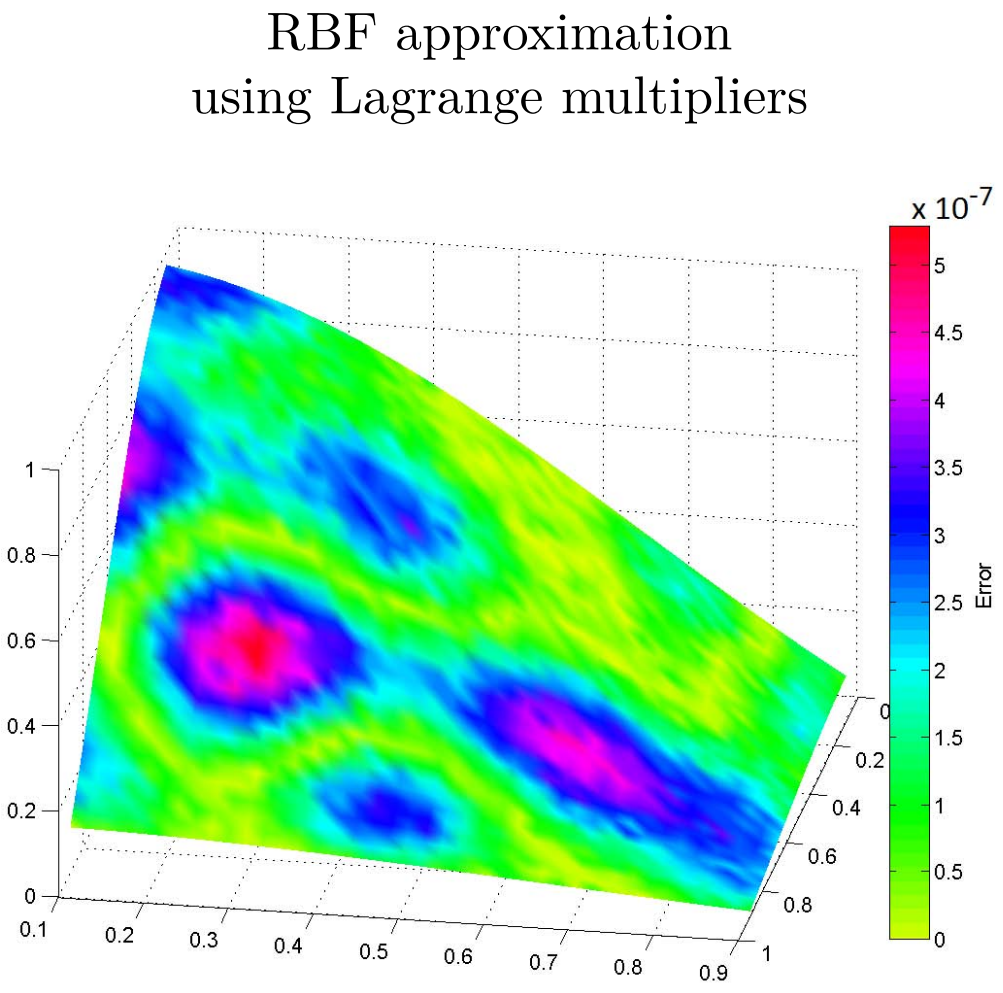

(c) Gauss, $\alpha = 1$, Halton points

Figure 5: Approximation to 1089 data points sampled from a $2D$ sinc function with 81 Halton-spaced Gaussian basis functions false-colored by magnitude of absolute error.

It can be seen that the RBF approximation using Lagrange multiplies (Fasshauer [10]) returns the worst result in terms of the error in comparison with the proposed methods. Further, in Fig. 6, it can be seen that the errors for all RBF approximation methods are much higher when the TPS is used.

There is a question of how the RBF approximation depends on the shape parameter $\alpha$. This is described in the following section.

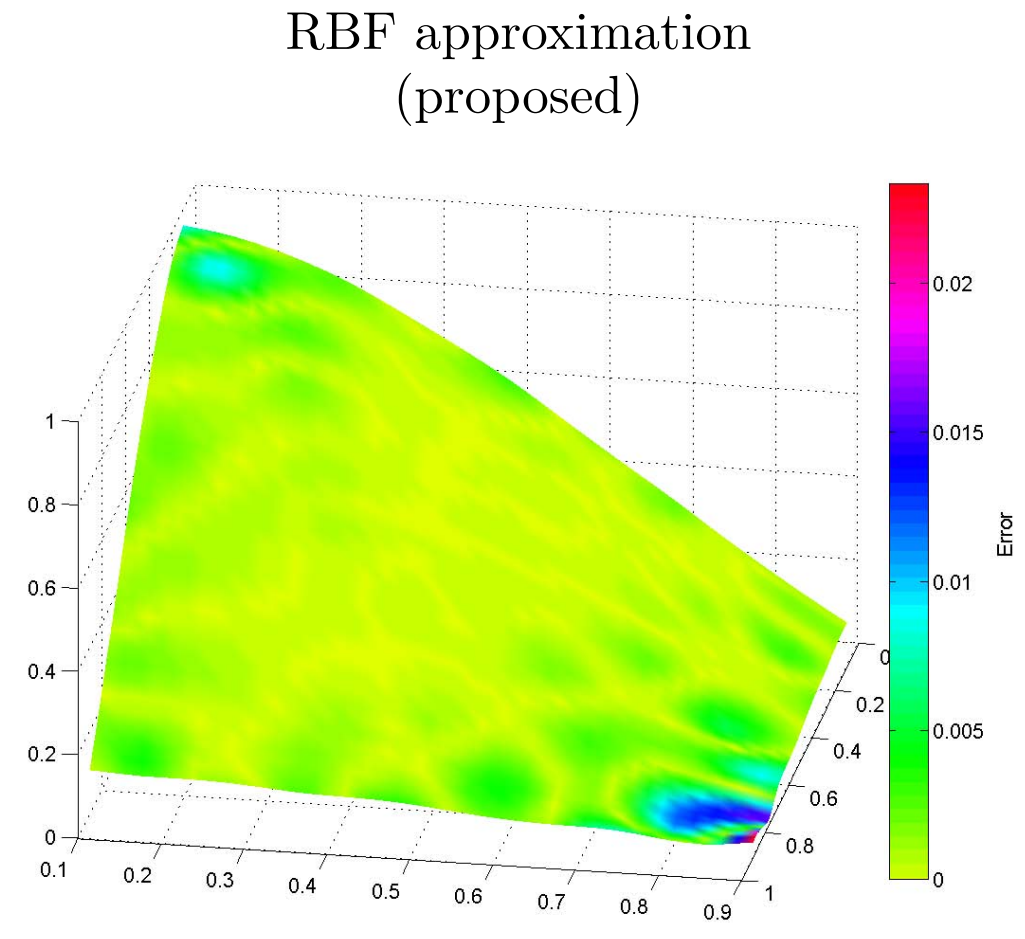

(a) TPS, $\alpha = 1$, Halton points

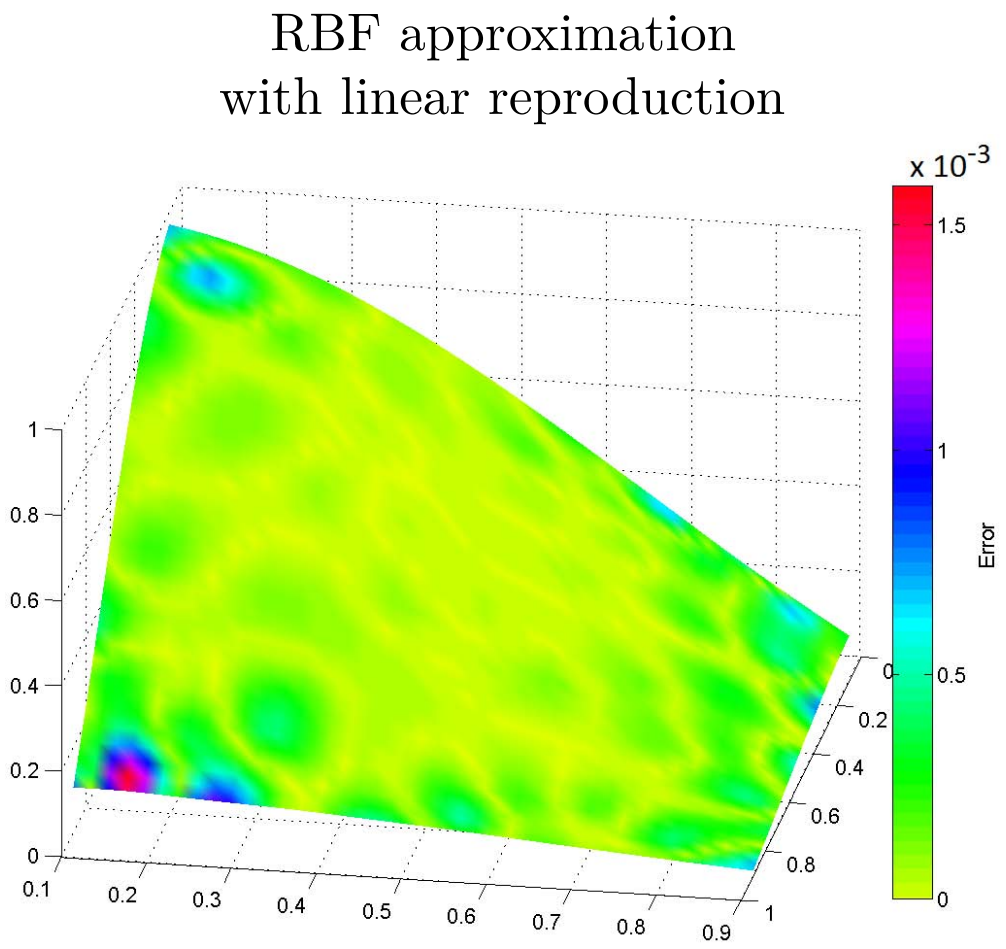

(b) TPS, $\alpha = 1$, Halton points

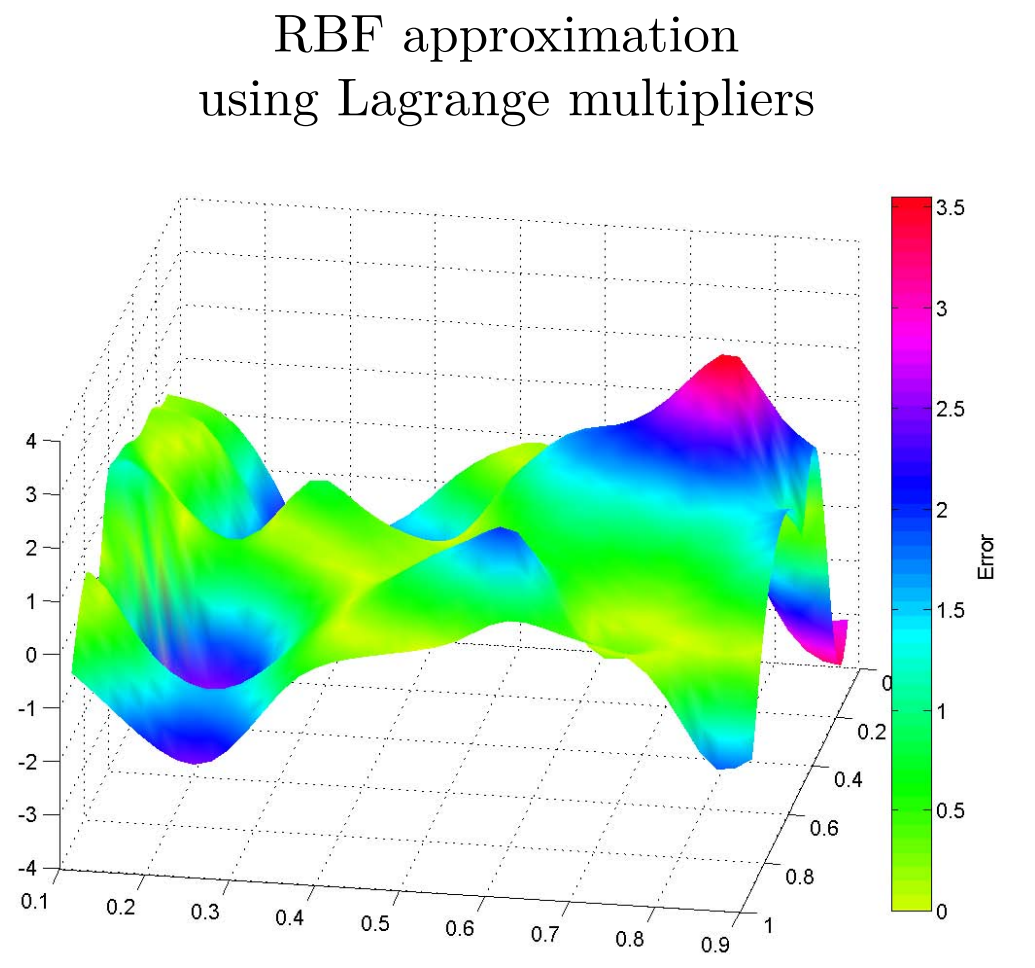

(c) TPS, $\alpha = 1$, Halton points

Figure 6: Approximation to 1089 data points sampled from a $2D$ sinc function with 81 Halton-spaced TPS false-colored by magnitude of absolute error.

### *5.2.2. Comparison of Methods*

In this section, the different versions of RBF approximation which were presented in Sect. 2 - Sect. 4 are compared. Figure 7 presents the mean absolute error of RBF approximation for the dataset, which consists of 1089 Halton points in the range $[0, 1] \times [0, 1]$, sampled from a $2D$ sinc function, while the set of reference points contains 81 points with Halton behavior of the distribution, and for different global radial basis functions. The graphs represent the mean absolute error according to a shape parameter $\alpha$ of used RBFs. We can see that for RBF approximation using Lagrange multipliers (Fasshauer [10]) we obtain a higher mean absolute error. Mean absolute errors for RBF approximation and RBF approximation with linear reproduction are almost the same. Moreover, the Gaussian RBF gives the best result for shape parameter $\alpha = 1$ and the inverse quadric for $\alpha = 0.5$. Further, the TPS function is not appropriate to solve the given problem, see Fig. 7c. Note, the standard deviation of errors was also measured and the same behavior and order of magnitude was obtained as for the mean absolute errors.

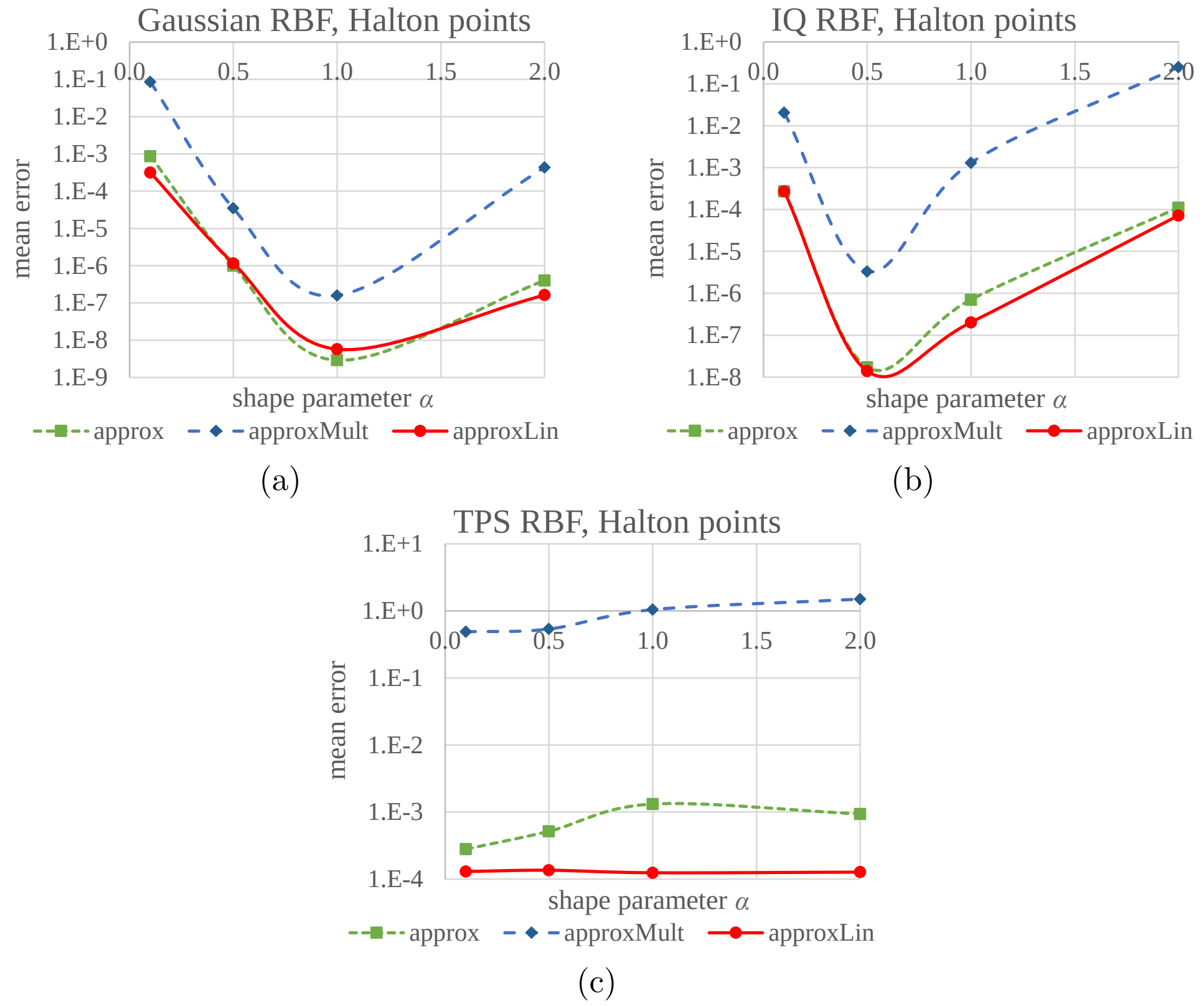

Figure 7: The mean absolute error of approximation to 1089 data points sampled from a $2D$ sinc function with 81 reference Halton points for different RBF approximation methods, different RBFs and different shape parameters. The used approximation methods are: proposed RBF approximation (approx), RBF approximation using Lagrange multipliers (approxMult) and RBF approximation with linear reproduction (approxLin). RBFs are: (a) Gauss function, (b) IQ, (c) TPS.

### *5.2.3. Comparison of Different Distributions of Reference Points*

In this section, we focus on a comparison of the presented RBF approximation methods due to used distribution of reference points. Measurements of errors were performed for different type of RBFs with different shape parameters. Mean absolute error according to shape parameter $\alpha$ for Gaussian RBF is presented in Fig. 8.

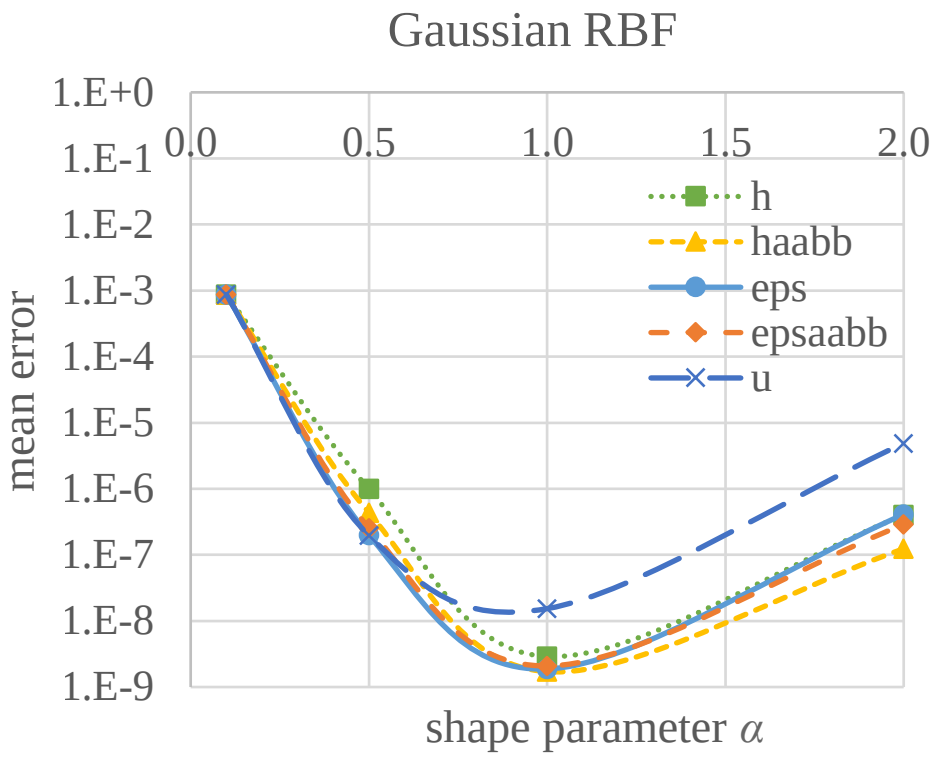

(a)

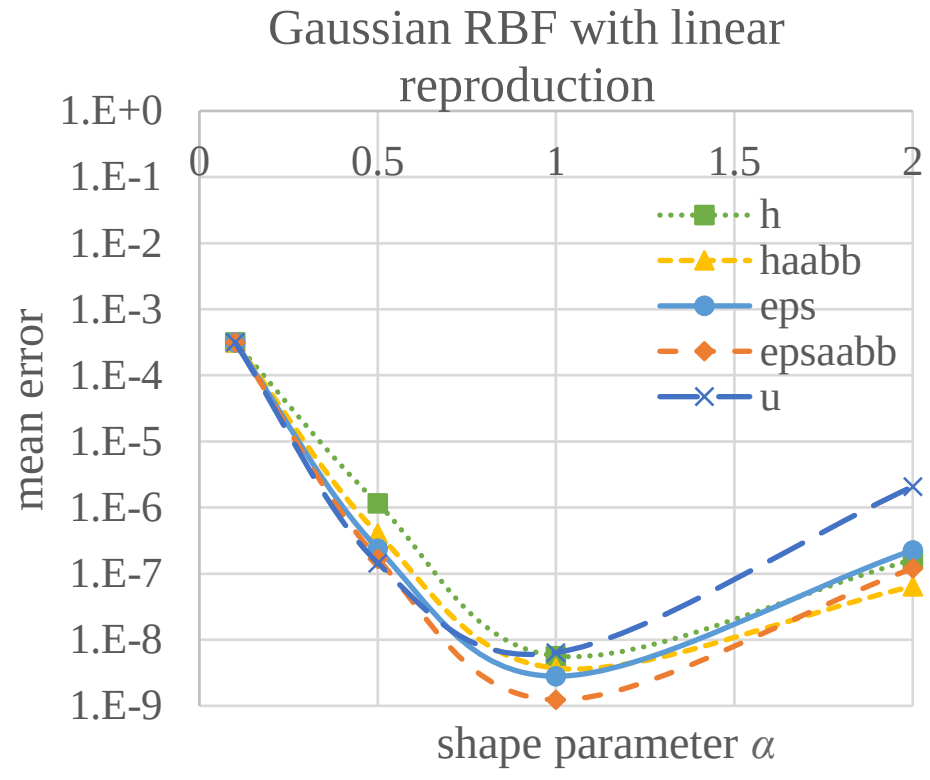

(b)

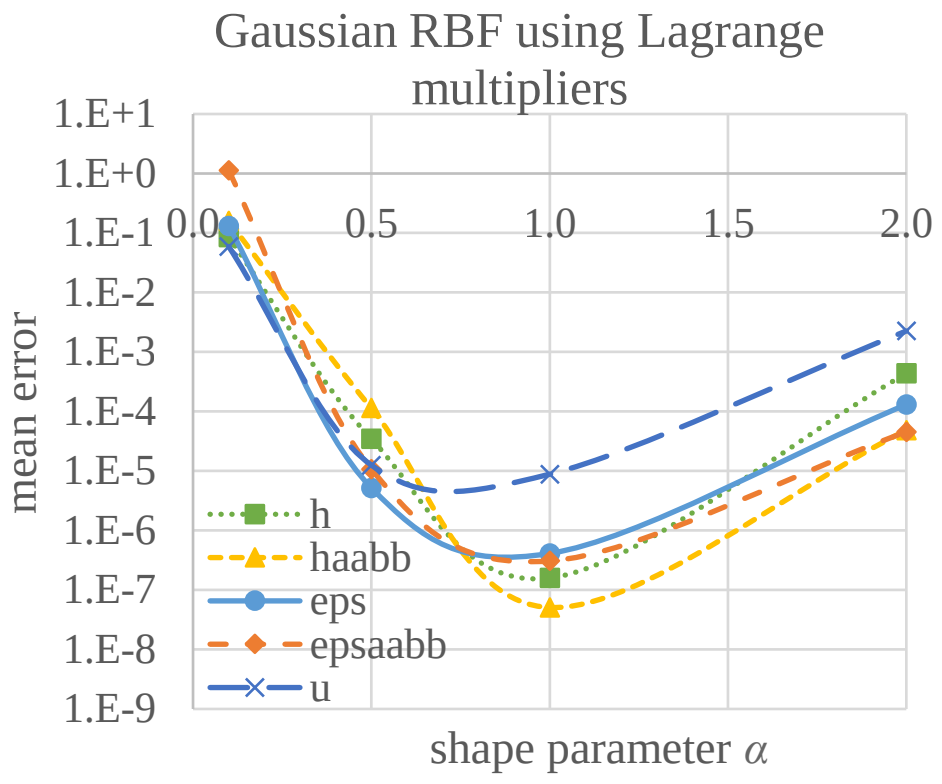

(c)

Figure 8: The mean absolute error of approximation to 1089 data points sampled from a $2D$ sinc function with 81 spaced Gaussian basis functions for different RBF approximation methods, different shape parameters and different sets of reference points. The sets of reference points are: Halton points (h), Halton points + AABB (haabb), epsilon points (eps), epsilon points + AABB (epsaabb), points on a regular grid (u). Their description is in Sect. 5.1. Versions of approximation are: (a) RBF approximation, (b) RBF approximation with linear reproduction, (c) RBF approximation using Lagrange multipliers.

We can see that for all versions of RBF approximation the worst result is obtained for reference points on a regular grid (u). For the proposed RBF approximation, the remaining sets of reference points give almost the same results. Reference points corresponding to epsilon points + AABB (epsaabb) almost always give the best result for RBF approximation with linear reproduction. For RBF approximation using Lagrange multipliers, the best results are for the reference points which have a Halton distribution.

### *5.2.4. Comparison by Placement of the Dataset in $E^2$*

This section is focused on placement of the actual dataset in the domain space and the used function generating associated scalar values in $E^2$. The given dataset has a range of one in both axes and the function generating associated scalar values is a $2D$ sinc function. Two configurations for placement of the origin of the dataset and the maximum of the $2D$ sinc function were used. The first configuration is at point $(0; 0)$; the second is moved to point $(3, 951, 753; 2, 785, 412)$.

Figure 9 presents the mean absolute error for these configurations, when the Gaussian basis functions and Halton set of reference points were chosen. We can see that RBF approximation with linear reproduction gives a higher error for the second configuration, i.e. placement at point $(3, 951, 753; 2, 785, 412)$.

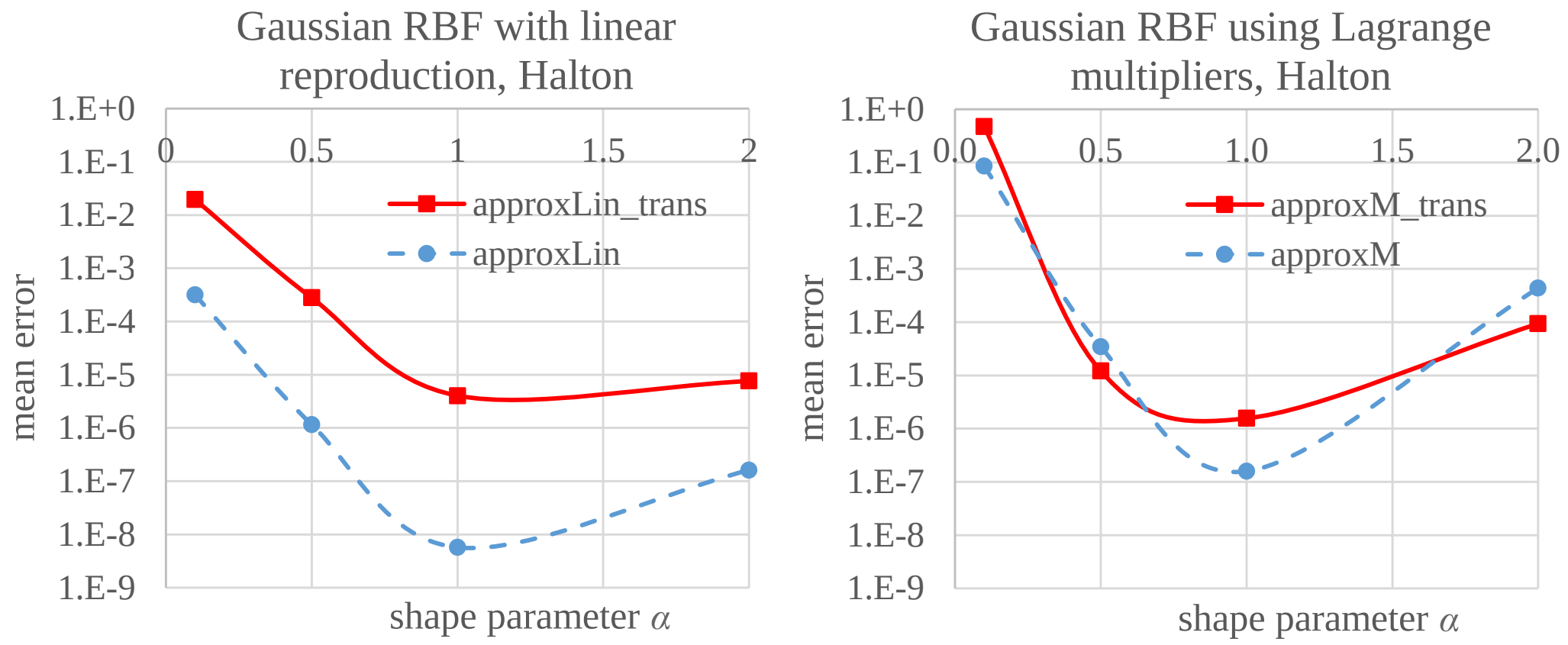

Figure 9: The mean absolute error of approximation to 1089 data points sampled from a $2D$ sinc function with 81 spaced Gaussian basis functions for a Halton set of reference points, different RBF approximation methods and different shape parameters. The placement of the given dataset and the maximum of the $2D$ sinc function are at point $(0; 0)$ (circles) or at point $(3, 951, 753; 2, 785, 412)$ (squares). Versions of approximation are RBF approximation with linear reproduction (left) and RBF approximation using Lagrange multipliers (right).

For RBF approximation using Lagrange multiplier the decision is not ambiguous. Note that a graph for the proposed RBF approximation is not presented, because both configurations give the same results.

*5.2.5. Optimal Number of Reference Points*

This section focuses on the influence of the number of reference points. The number of reference points is determined relative to the number of points in the given dataset. Measurements for different shape parameters were performed many times and average mean absolute errors were computed, see Fig. 10 - Fig. 12. Note that the reference points were distributed by Halton distribution. Figure 10 presents the mean absolute error for the Gaussian RBF approximation. Experimental results for the IQ are shown in Fig. 11. We can see that for the small shape parameter $\alpha$ the mean absolute errors are almost constant. However, for greater shape parameters the mean absolute error decreases with the increasing number of reference points.

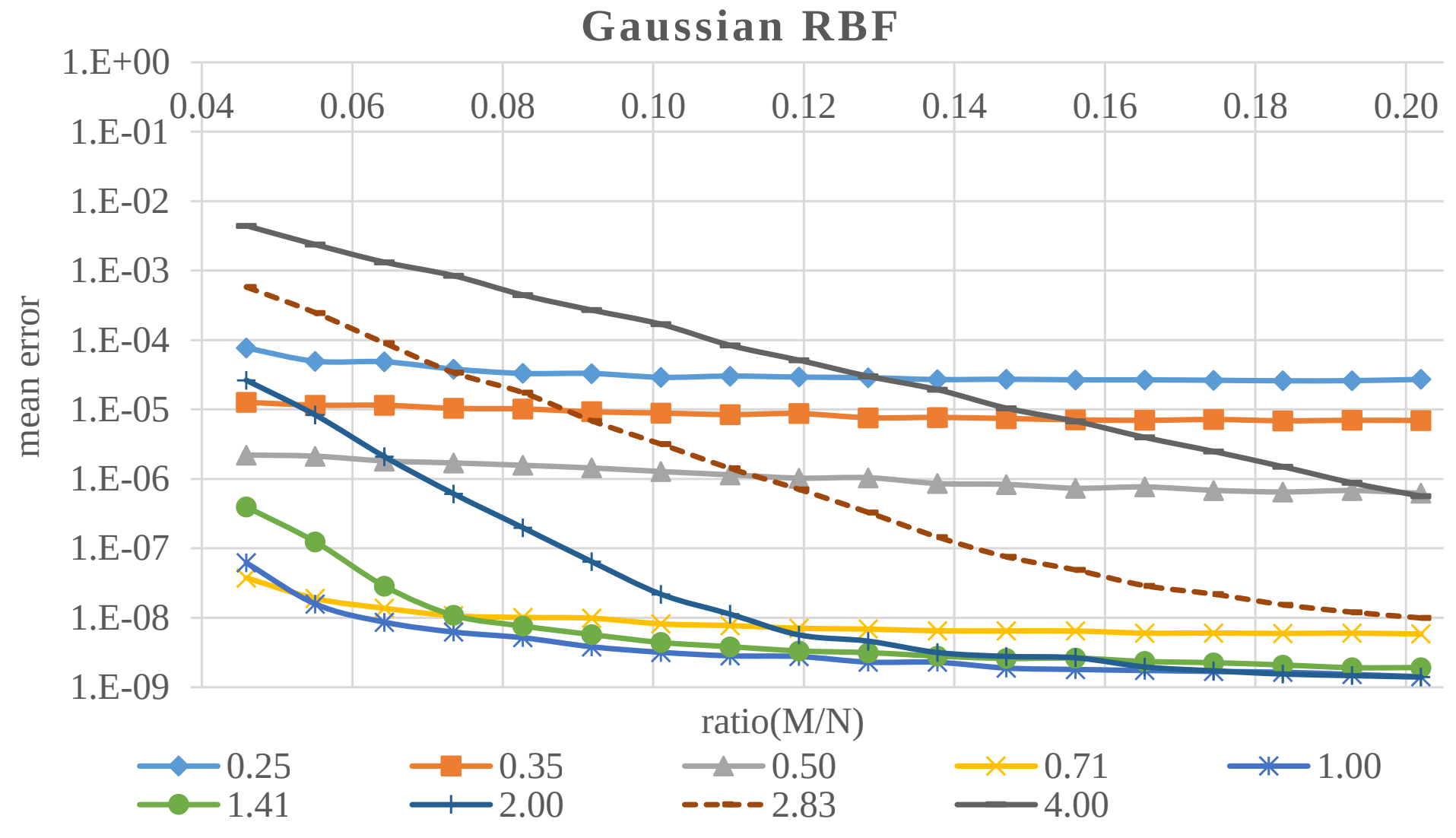

Figure 10: The mean absolute error of the proposed RBF approximation to 1089 data points sampled from a $2D$ sinc function for different numbers of reference points, Gaussian RBF with different shape parameters $\alpha$.

Figure 12 presents experimental results obtained for the TPS function. We can see that the mean absolute error decreases with the increasing number of reference points as would be expected.

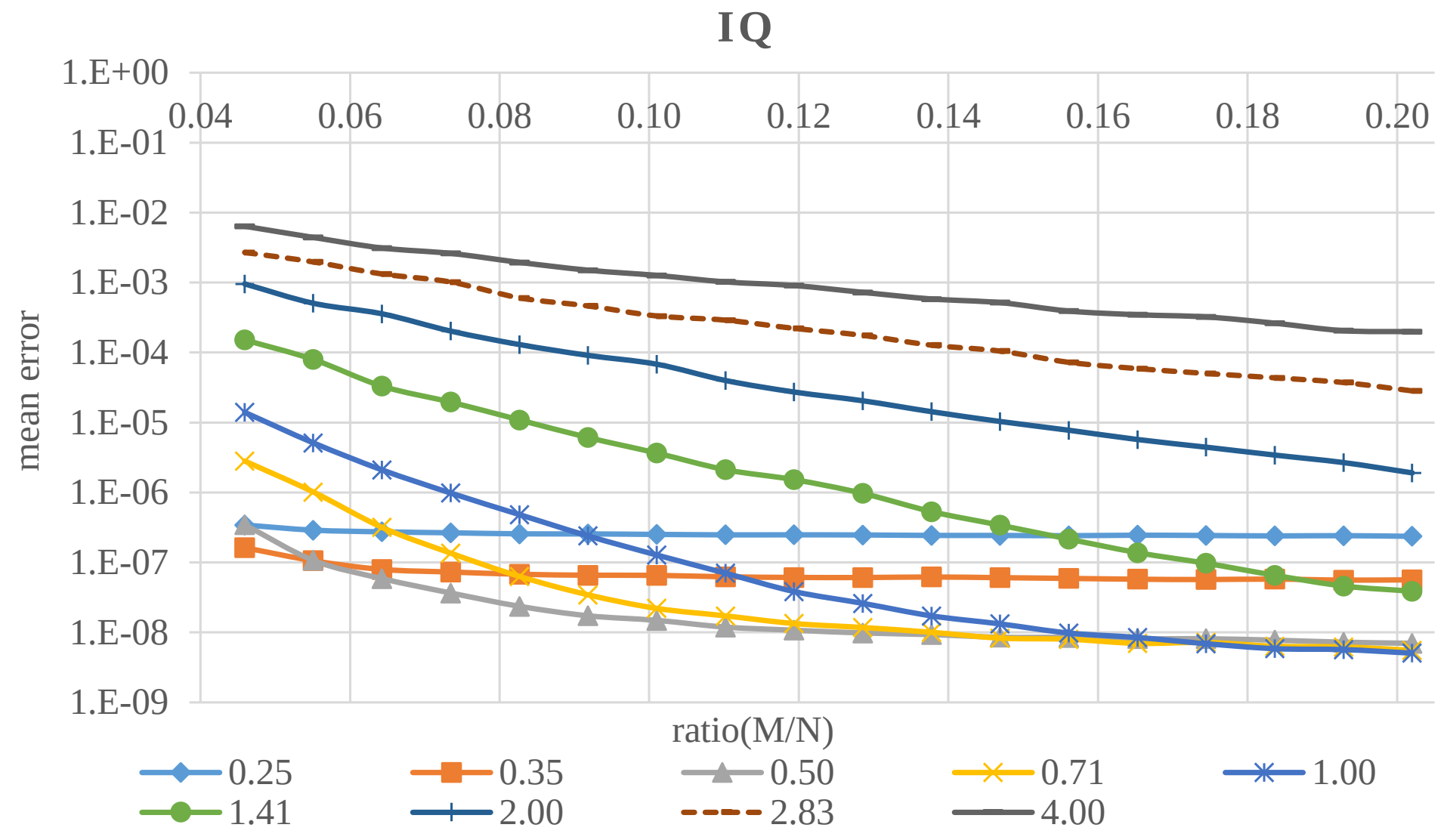

Figure 11: The mean absolute error of the proposed RBF approximation to 1089 data points sampled from a $2D$ sinc function for different numbers of reference points, IQ RBF with different shape parameters $\alpha$.

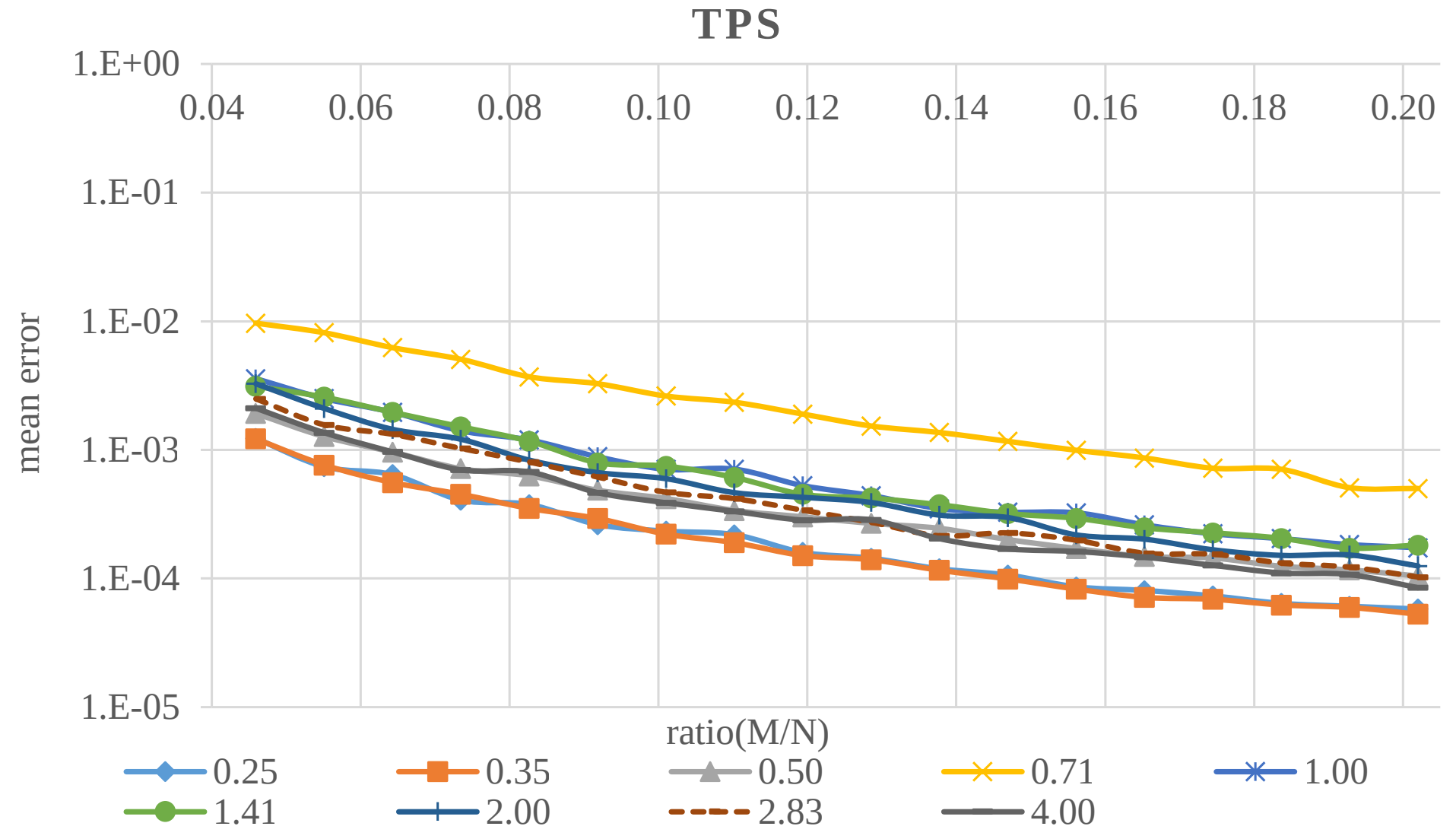

Figure 12: The mean absolute error of the proposed RBF approximation to 1089 data points sampled from a $2D$ sinc function for different numbers of reference points, TPS RBF with different shape parameters $\alpha$.

Finally, note that the results for RBF approximation with reproduction are very similar to the proposed RBF approximation. RBF using Lagrange multipliers has unpredictable behavior and no trend can be established.

### *5.3. Real Datasets*

The presented methods of the RBF approximation have been also tested on real data. Let us introduce results for real dataset which was obtained from GPS data of mount Veľký Rozsutec in the Malá Fatra, Slovakia[1]. Each point of this dataset is associated with its elevation. Moreover, as a first step, the real dataset is translated so that its estimated center of gravity corresponds to the origin of the coordinate system. This step is used due to the limitation of the influence of dataset placement in space and it was chosen based on the results of experiments described in Sect. 5.2.4. Table 2 gives an overview of the used dataset.

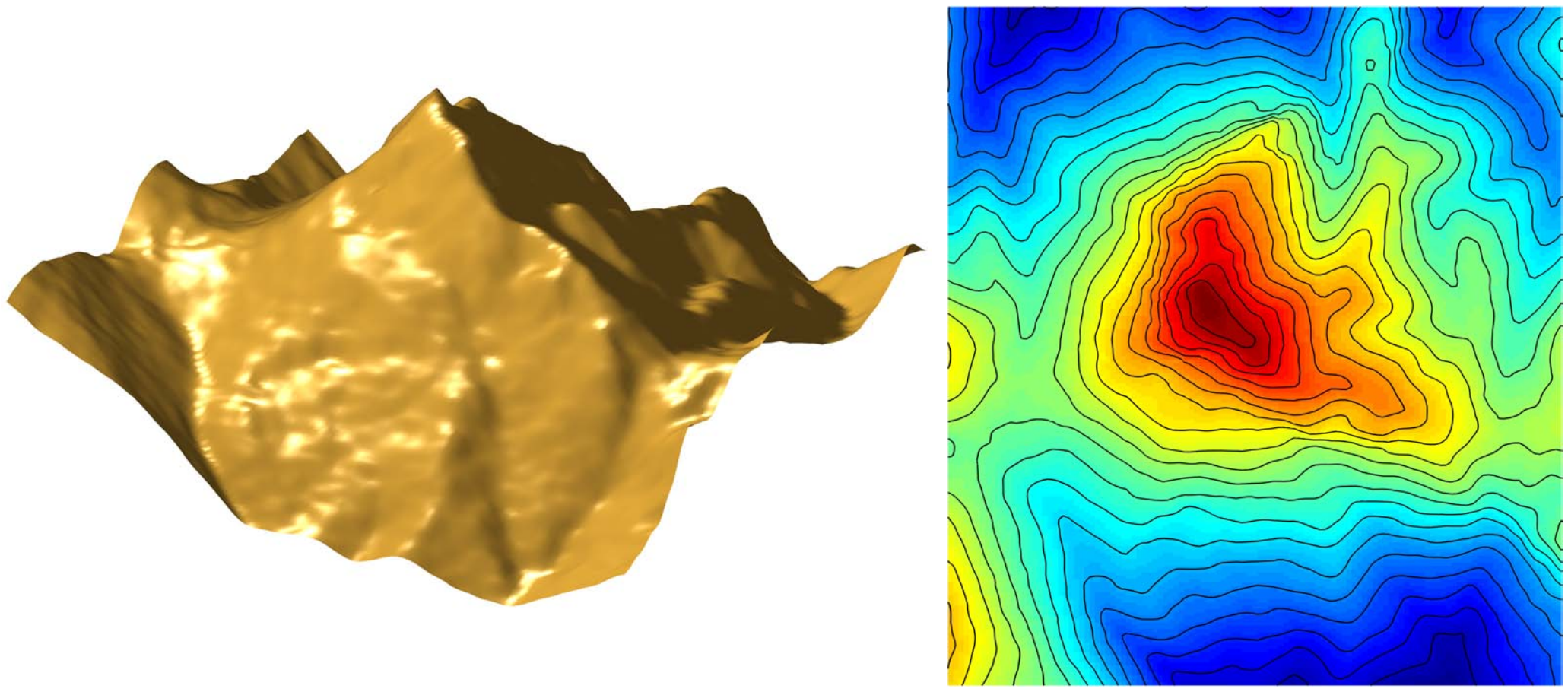

Figure 13: Mount Veľký Rozsutec, Slovakia (left) and its contour map (right).

Table 2: Overview information for the tested real dataset. The Axis-Aligned Bounding Box (AABB) of the tested dataset has a size width × length × relief, i.e. $x_{range} \times y_{range} \times z_{range}$.

| | **Veľký Rozsutec** |
|---|---|
| **number of pts.** | 24, 190 |
| **relief [m]** | 818.8000 |
| **width [m]** | 2608.5927 |
| **length [m]** | 2884.1169 |

[1]http://www.gpsvisualizer.com/elevation

#### *5.3.1. Examples of RBF Approximation Results*

Results for RBF approximation of mount Veľký Rozsutec dataset using Halton set of reference points, which consists of 484 points, and Gaussian RBF with shape parameter $\alpha = 0.0025$ are shown in Fig. 14 and histograms of errors for these results are shown in Fig. 15.

Note, that the results of RBF approximation using Lagrange multipliers

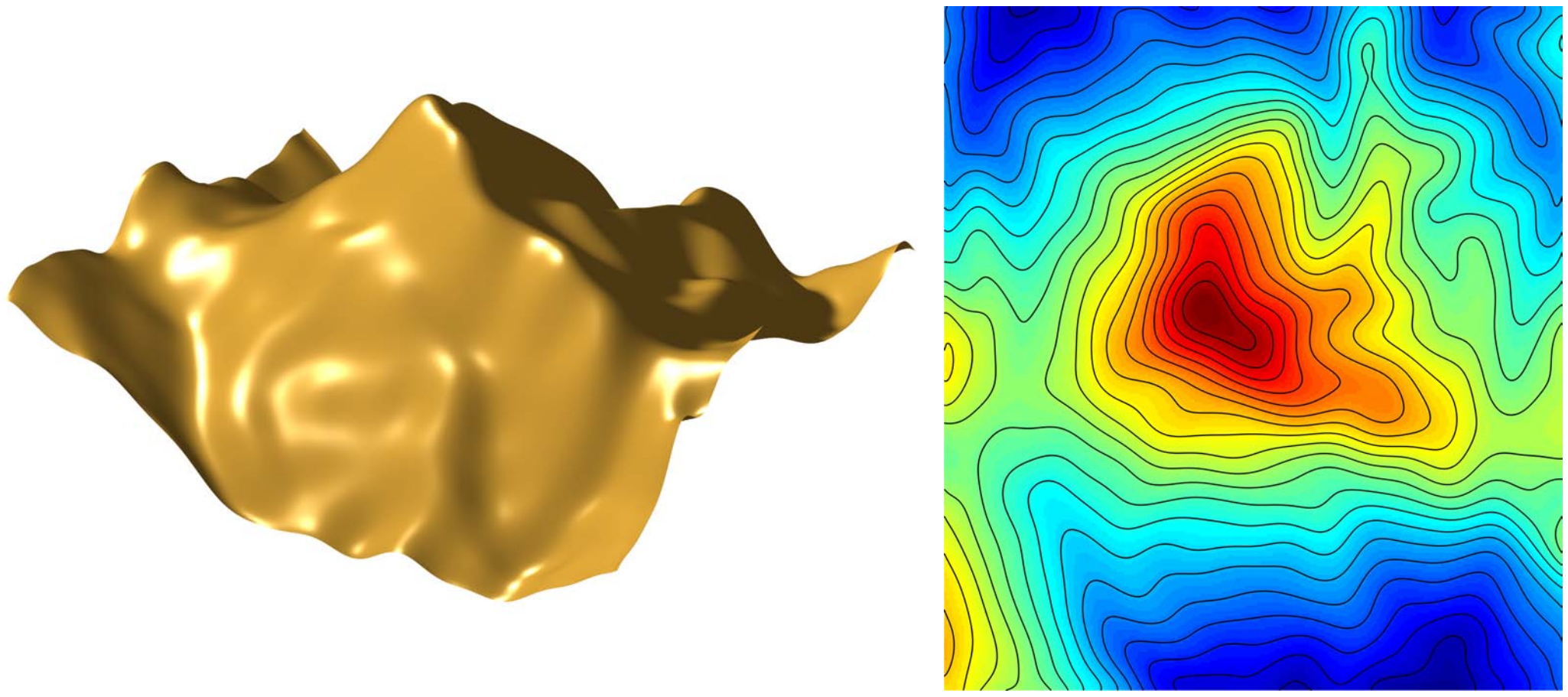

(a) RBF approximation: Gauss, $N = 24,190$, $M = 484$, $\alpha = 0.0025$, Halton points (left) and its contour map (right)

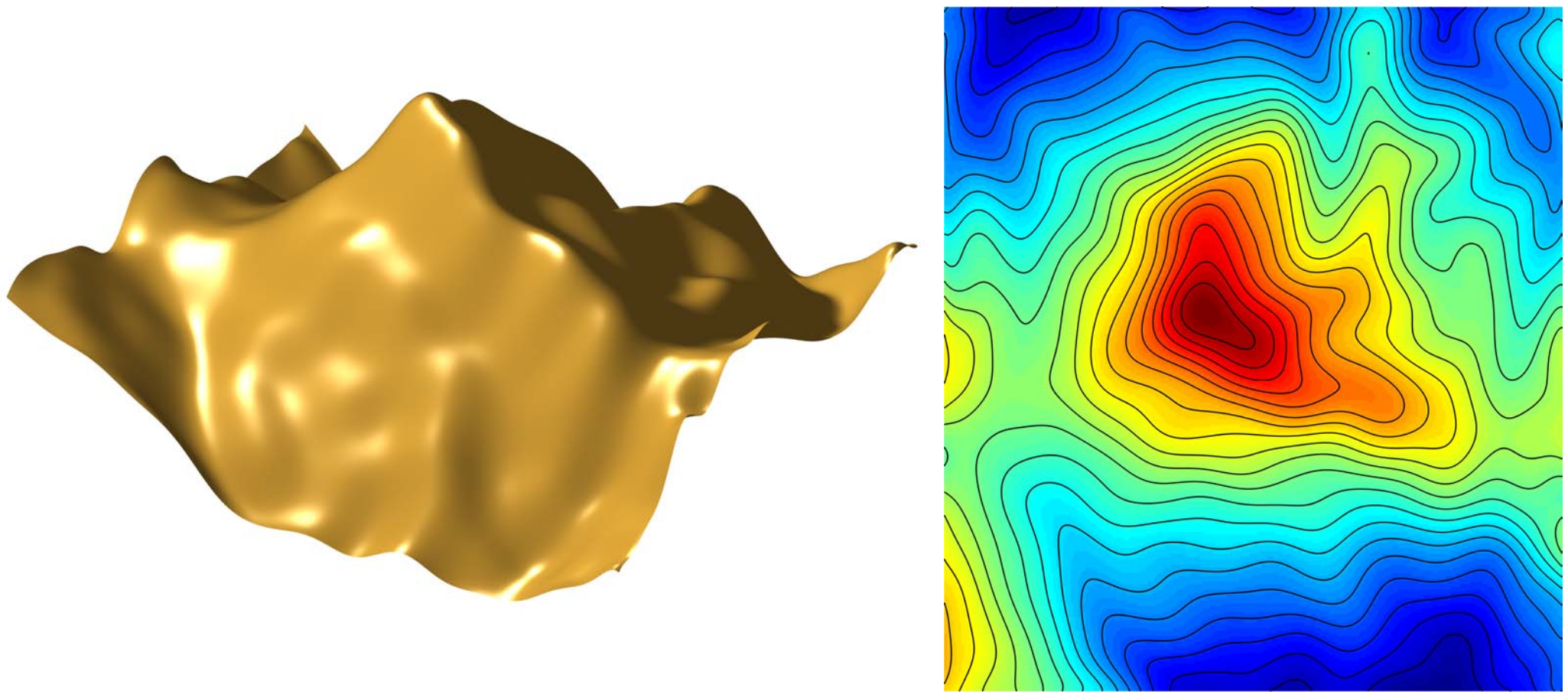

(b) RBF approximation with reproduction: Gauss, $N = 24,190$, $M = 484$, $\alpha = 0.0025$, Halton points (left) and its contour map (right) are not presented for real data because this method has unpredictable behavior and is unusable for real dataset, which was already evident from results for synthetic datasets. From presented results, it can be seen that the RBF approximation with linear reproduction returns the worst result in terms of the error in comparison with the proposed method. Moreover, if the results of approximation are compared with the original, it can be seen that the RBF approximation with the global Gaussian RBFs cannot preserve the sharp ridge.

Figure 14: Results for mount Veľký Rozsutec approximated by 484 Halton-spaced Gaussian basis functions with shape parameter $\alpha = 0.0025$.

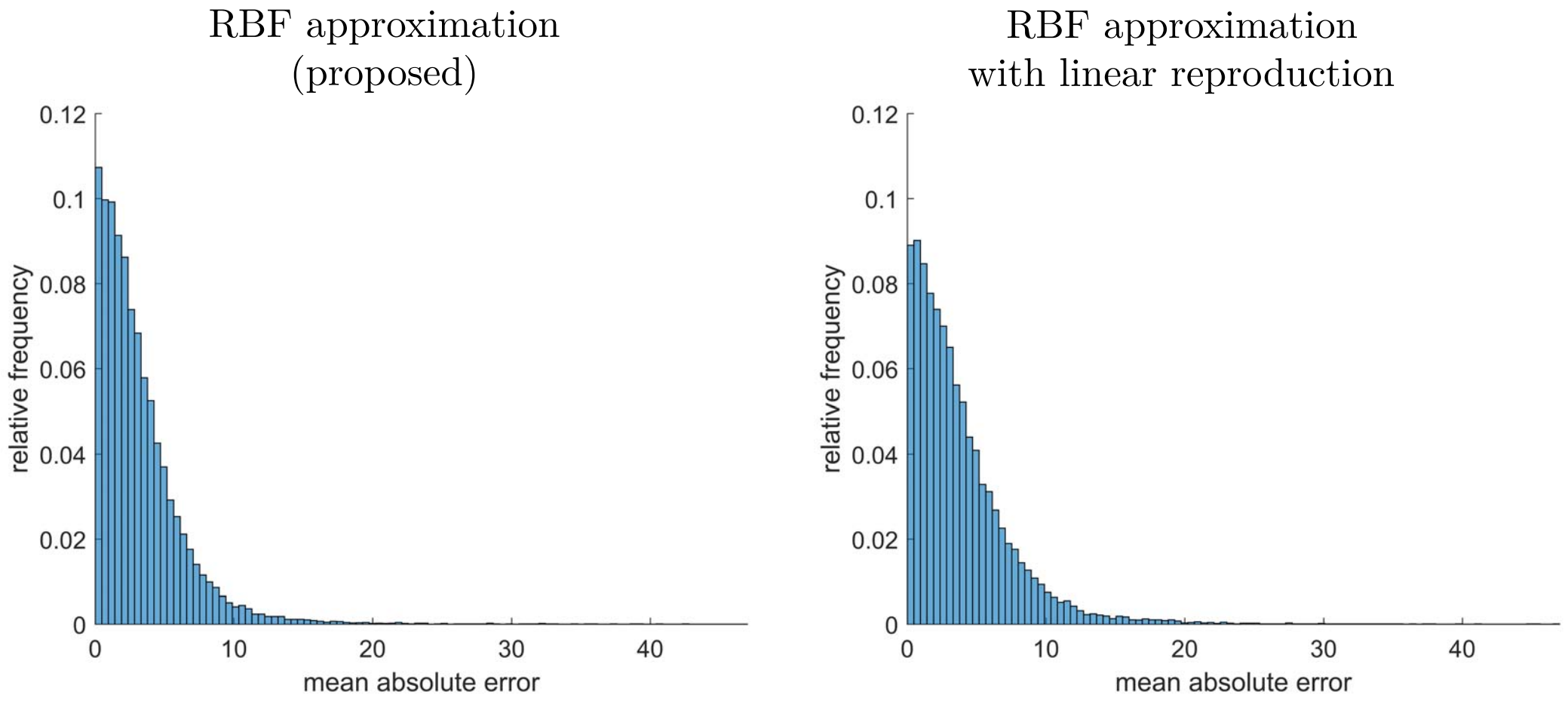

Figure 15: Histograms of errors for mount Veľký Rozsutec approximated by 484 Halton-spaced Gaussian basis functions with shape parameter $\alpha = 0.0025$.

Results for RBF approximation of mount Veľký Rozsutec dataset using Halton set of reference points, which contains different number of points, and TPS with shape parameter $\alpha = 0.005$ are show in Fig. 16. The histograms of errors for these results are shown in Fig. 17. From these results, it can be seen that with an increasing number of reference points, approximation error is improved and some surface details also begin to appear. However, it can be again seen that the RBF approximation with the global TPS cannot preserve the sharp ridge.

There is a question of how the RBF approximation of real dataset depends on the shape parameter $\alpha$ and distribution of reference points. This is described in the following sections.

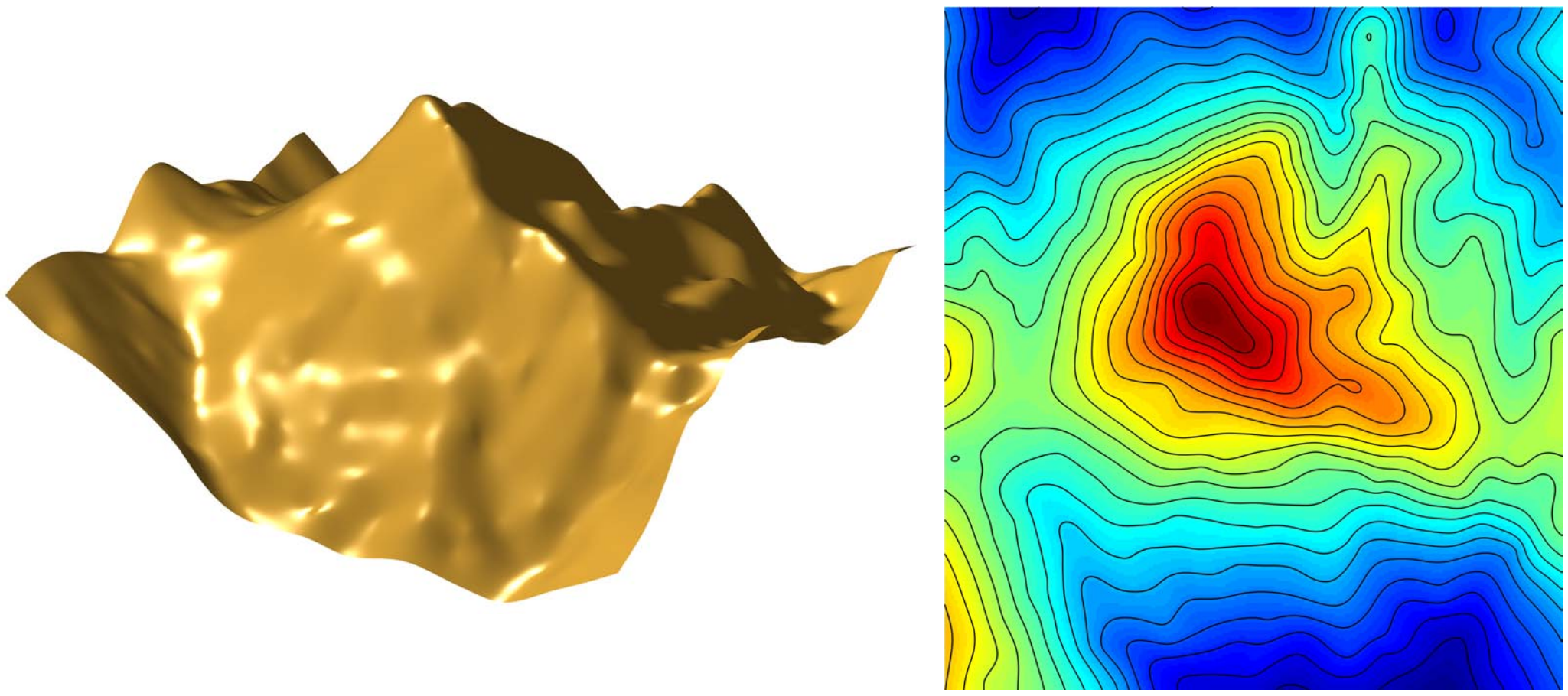

(a) RBF approximation: TPS, $N = 24,190$, $M = 484$, $\alpha = 0.005$, Halton points (left) and its contour map (right)

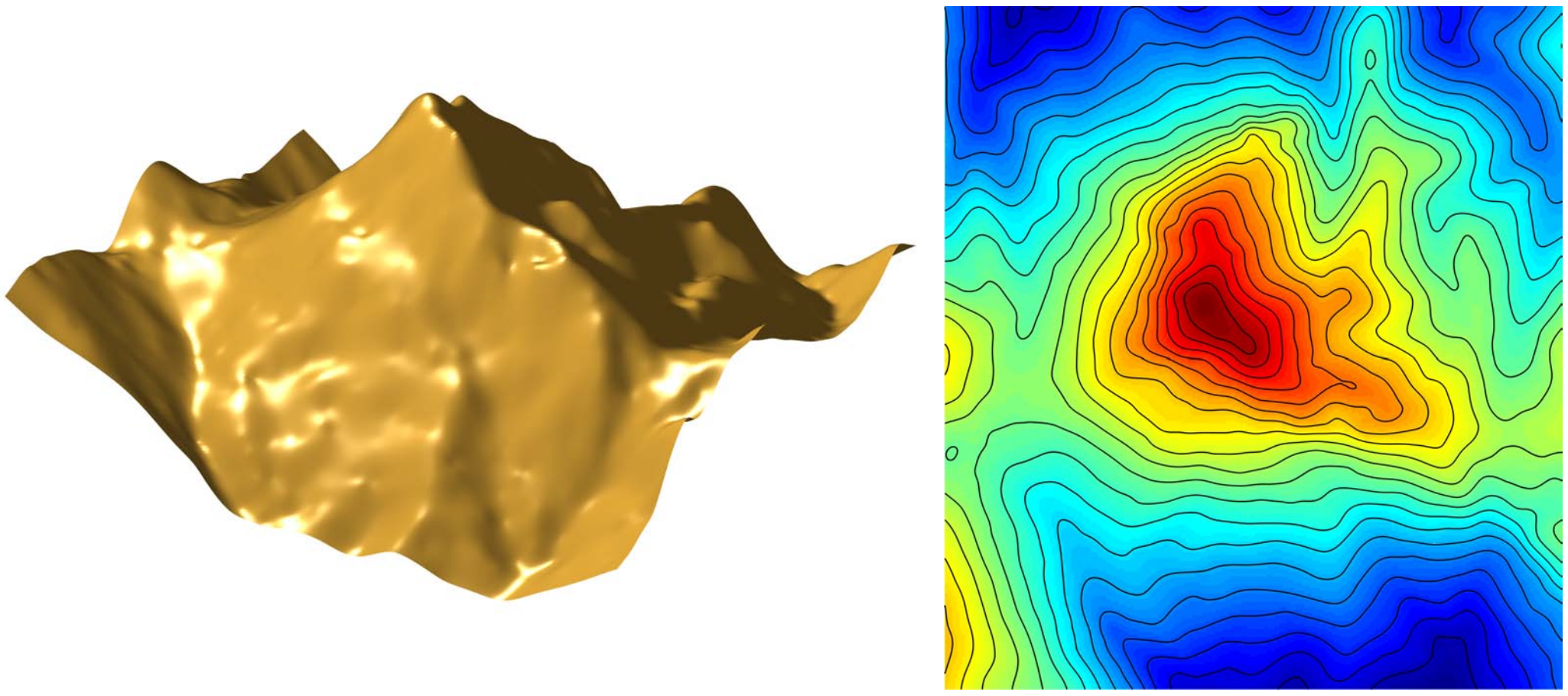

(b) RBF approximation: TPS, $N = 24,190$, $M = 1089$, $\alpha = 0.005$, Halton points (left) and its contour map (right)

Figure 16: Results for mount Veľký Rozsutec approximated by Halton-spaced TPS with shape parameter $\alpha = 0.005$.

### *5.3.2. Comparison of Different Distributions of Reference Points*

In this section, we focus on a comparison of the presented RBF approximation methods due to used distribution of reference points when the real data are approximated. Measurements of errors were performed for different type of RBFs with different shape parameters. Mean relative error according to shape parameter $\alpha$ for the Gaussian RBF is presented in Fig. 18 and for the IQ is shown in Fig. 19. Note that the mean *relative* error is presented for real data. The reason for this choice is that the function values of the real dataset are not normalized to the interval $[0, 1]$.

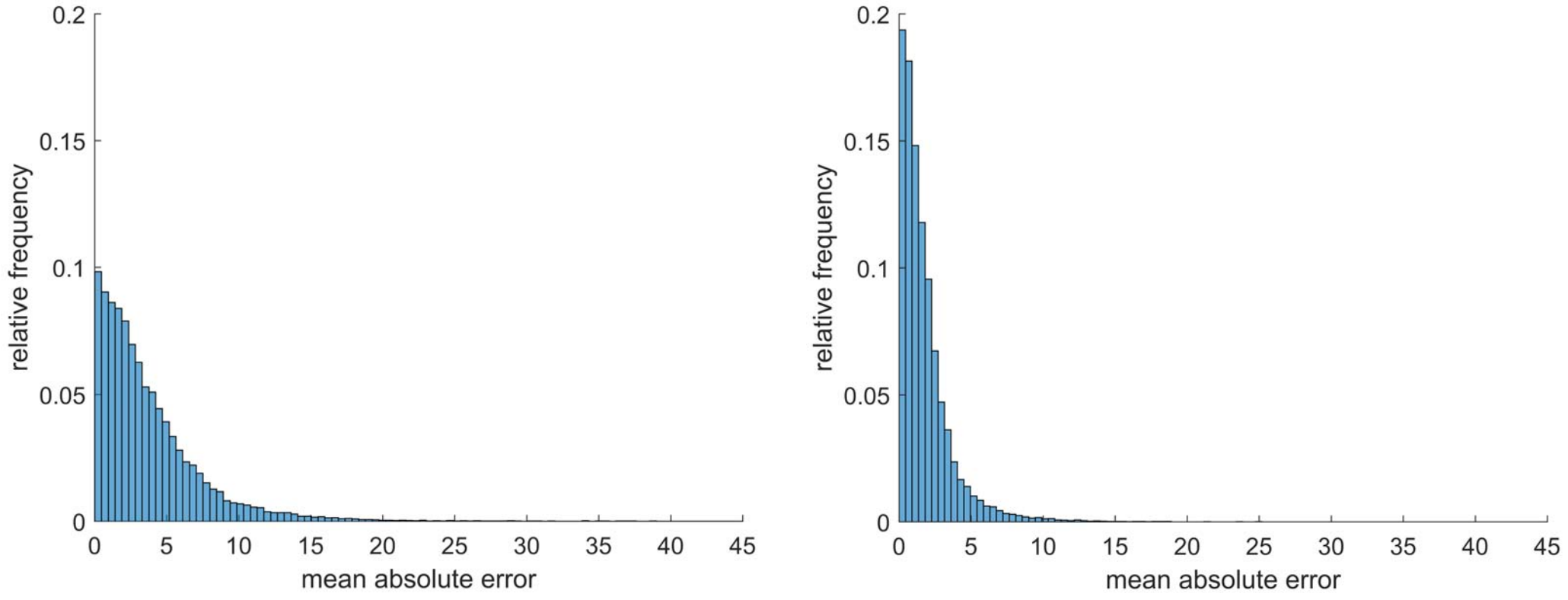

(a) $M = 484$ (b) $M = 1089$

Figure 17: Histograms of errors for mount Veľký Rozsutec approximated by Halton-spaced TPS with shape parameter $\alpha = 0.005$.

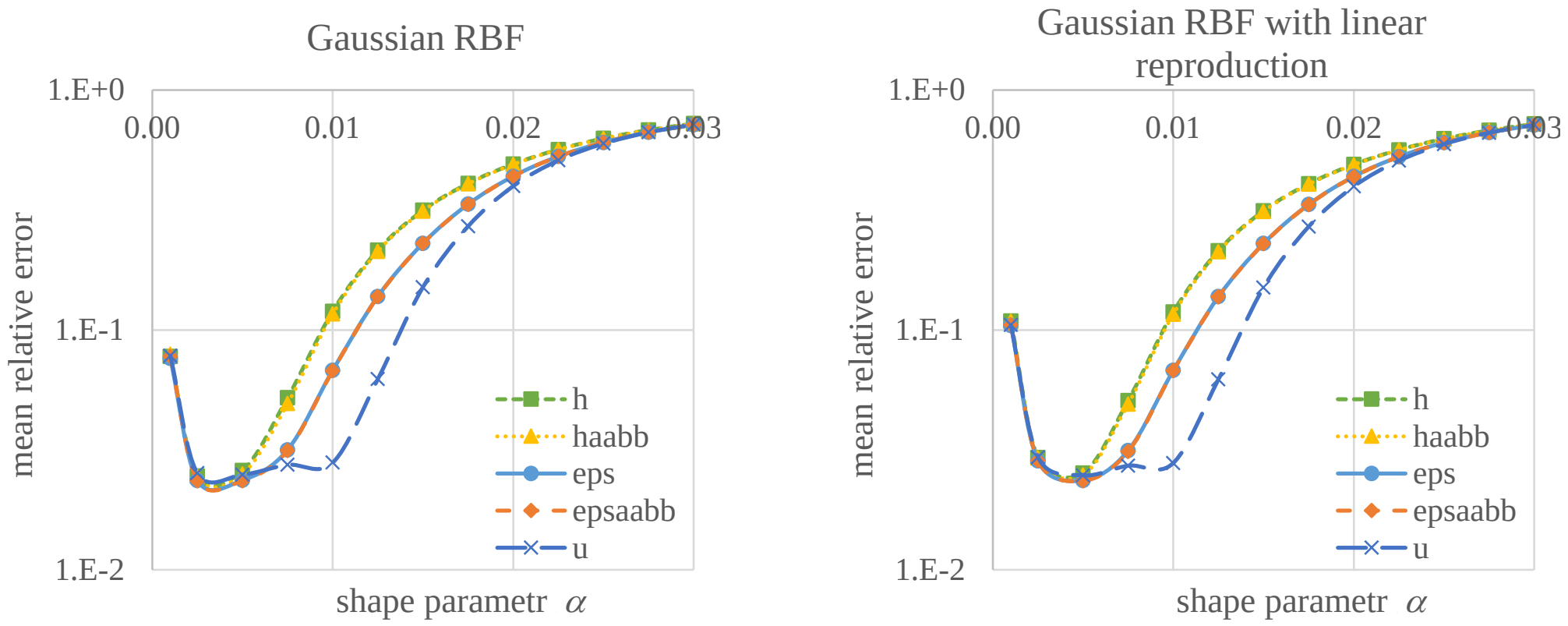

(a) without reproduction (b) with linear reproduction

Figure 18: The mean relative error of approximation for mount Veľký Rozsutec with 484 spaced Gaussian basis functions for different RBF approximation methods, different shape parameters and different sets of reference points. The sets of reference points are: Halton points (h), Halton points + AABB (haabb), epsilon points (eps), epsilon points + AABB (epsaabb), points on a regular grid (u), described in Sect. 5.1.

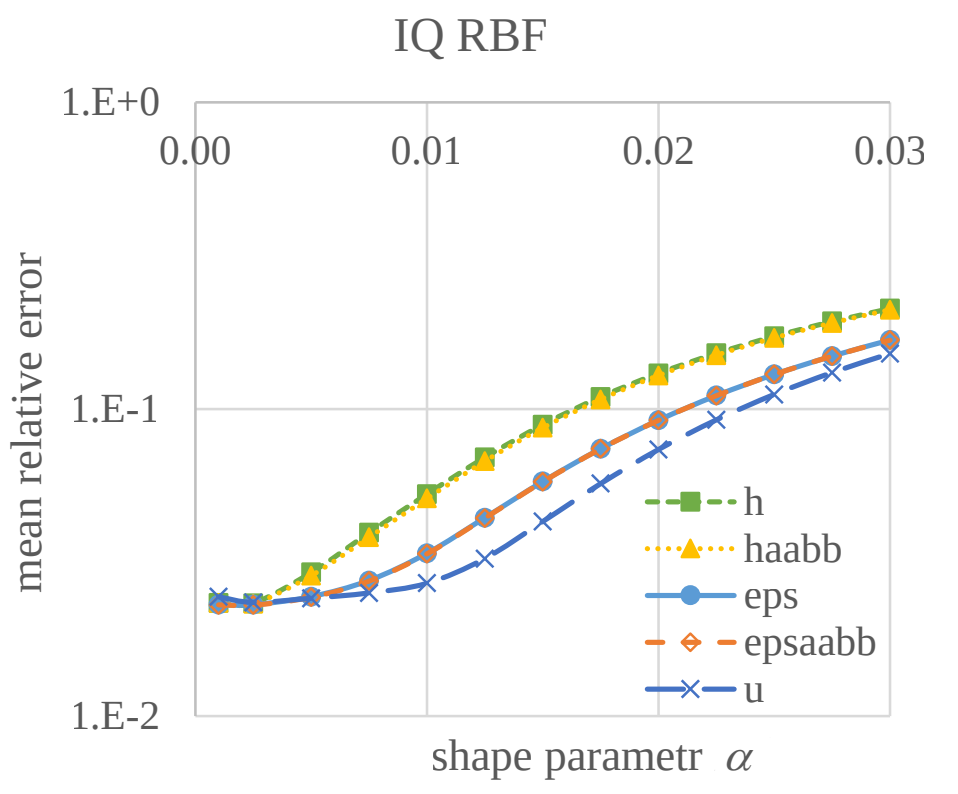

(a) without reproduction

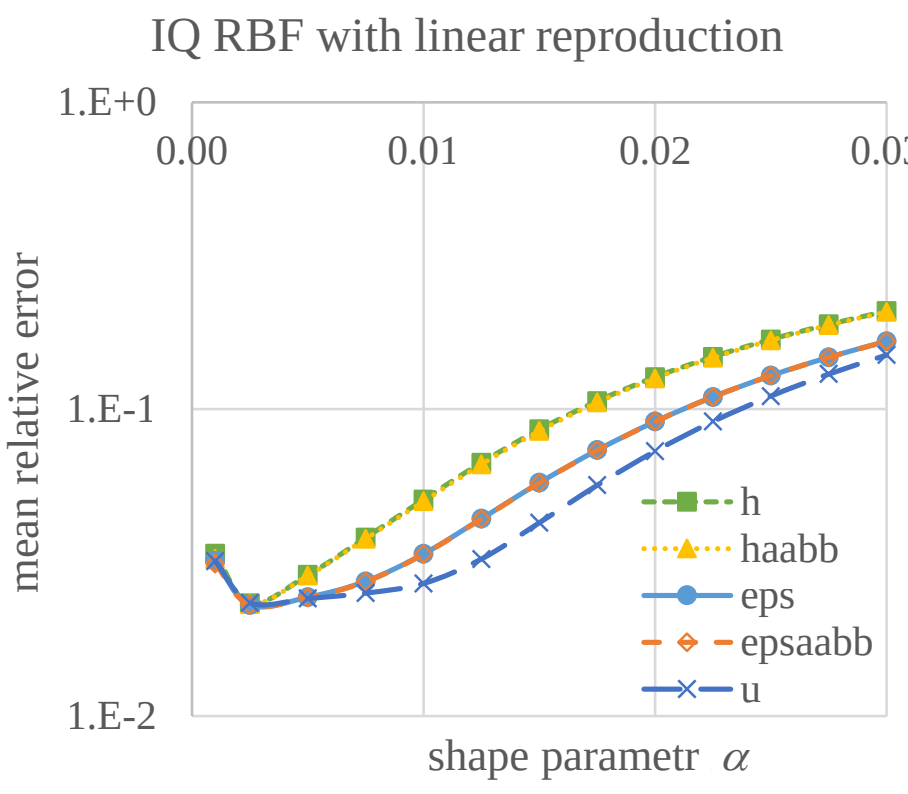

(b) with linear reproduction

Figure 19: The mean relative error of approximation for mount Veľký Rozsutec with 484 spaced IQ for different RBF approximation methods, different shape parameters and different sets of reference points. The sets of reference points are: Halton points (h), Halton points + AABB (haabb), epsilon points (eps), epsilon points + AABB (epsaabb), points on a regular grid (u), described in Sect. 5.1.

We can see that for all versions of RBF approximation, if the shape parameter $\alpha$ is not close to the optimum, the worst results are obtained for reference points with Halton distribution ((h) and (haabb)). The best results are obtained for reference points on a regular grid (u). If the shape parameter is chosen close to the optimum (for the presented configuration $\alpha \approx 0.0025$) then the mean relative error has only minor differences for different distribution of reference points. These results are different in comparison with results obtained for synthetic data.

Finally, note that the mean relative error for approximation of mount Veľký Rozsutec dataset according to shape parameter is constant for the TPS and deviation of mean relative error for different distribution of reference points is almost negligible.

### *5.3.3. Comparison of Different Radial Basis Functions*

In this section, we focus on a comparison of the results of RBF approximations using different types of RBFs. Real datasets were used for experiments and results of the mount Veľký Rozsutec are presented. Measurements of errors were performed for Halton set with 484 reference points. The shape parameter $\alpha = 0.0025$ was chosen for all types of RBFs. The differences of frequencies of errors are shown in Fig. 20. It can be seen that the best error is obtained for the RBF approximation using the IQ function. On the contrary, the worst error returns the RBF approximation using the TPS function.

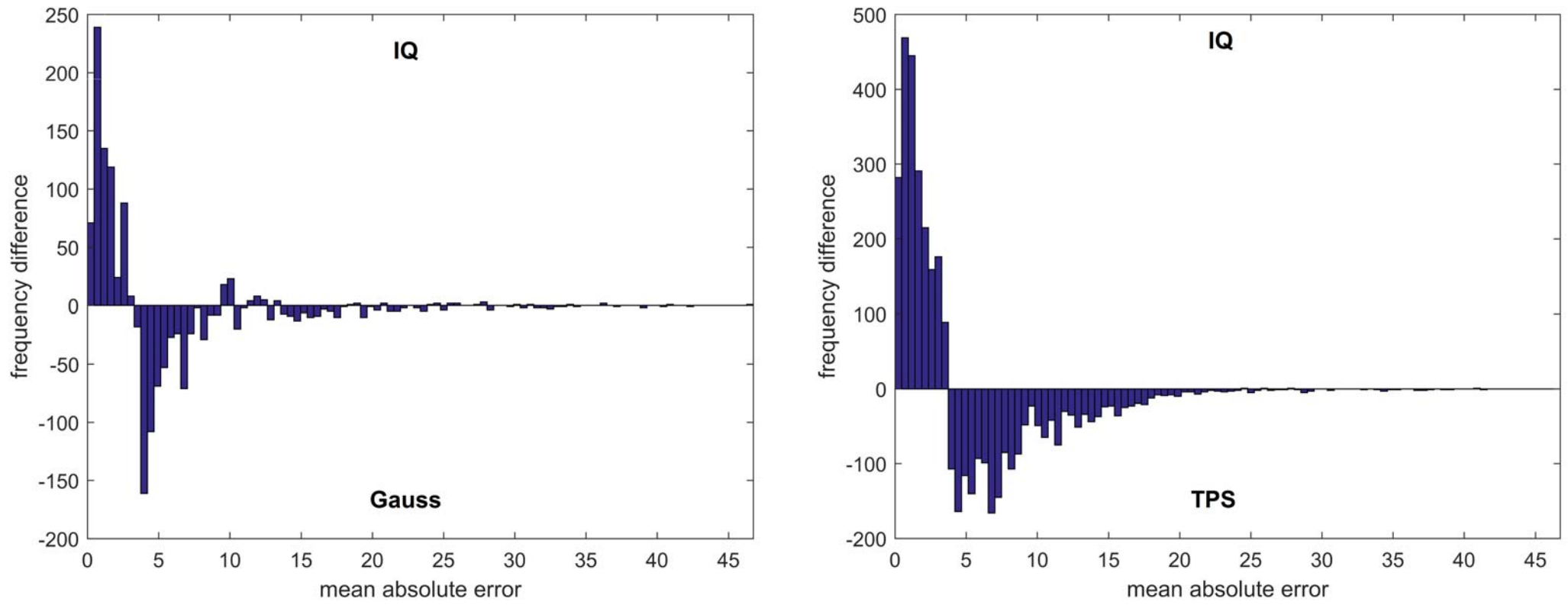

(a) IQ vs. Gaussian RBF (b) IQ vs. TPS

Figure 20: Difference of frequencies of error for mount Veľký Rozsutec approximated by 484 Halton-spaced RBFs with shape parameter $\alpha = 0.0025$.

Finally, note that the results for the RBF approximation with linear reproduction are similar to the proposed RBF approximation.

*5.3.4. Optimal Number of Reference Points*

This section focuses on the influence of the number of reference points for RBF approximation of mount Veľký Rozsutec dataset. The number of reference points is determined relative to the number of points in the given dataset. Measurements for different shape parameters were performed many times and average mean relative errors were computed, see Fig. 21 - Fig. 23. Note that the reference points were distributed by Halton distribution. Figure 21 presents the mean relative error for the Gaussian RBF approximation and Fig. 22 presents results for the IQ. It can be seen that for small shape parameter $\alpha$ the mean relative errors are almost constant. However, for the greater shape parameters the mean relative error decreases with the increasing number of reference points. These results are consistent with results for synthetic dataset.

Figure 23 presents experimental results obtained for TPS. We can see that the mean relative error is independent on the shape parameter $\alpha$ and decreases with the increasing number of reference points.

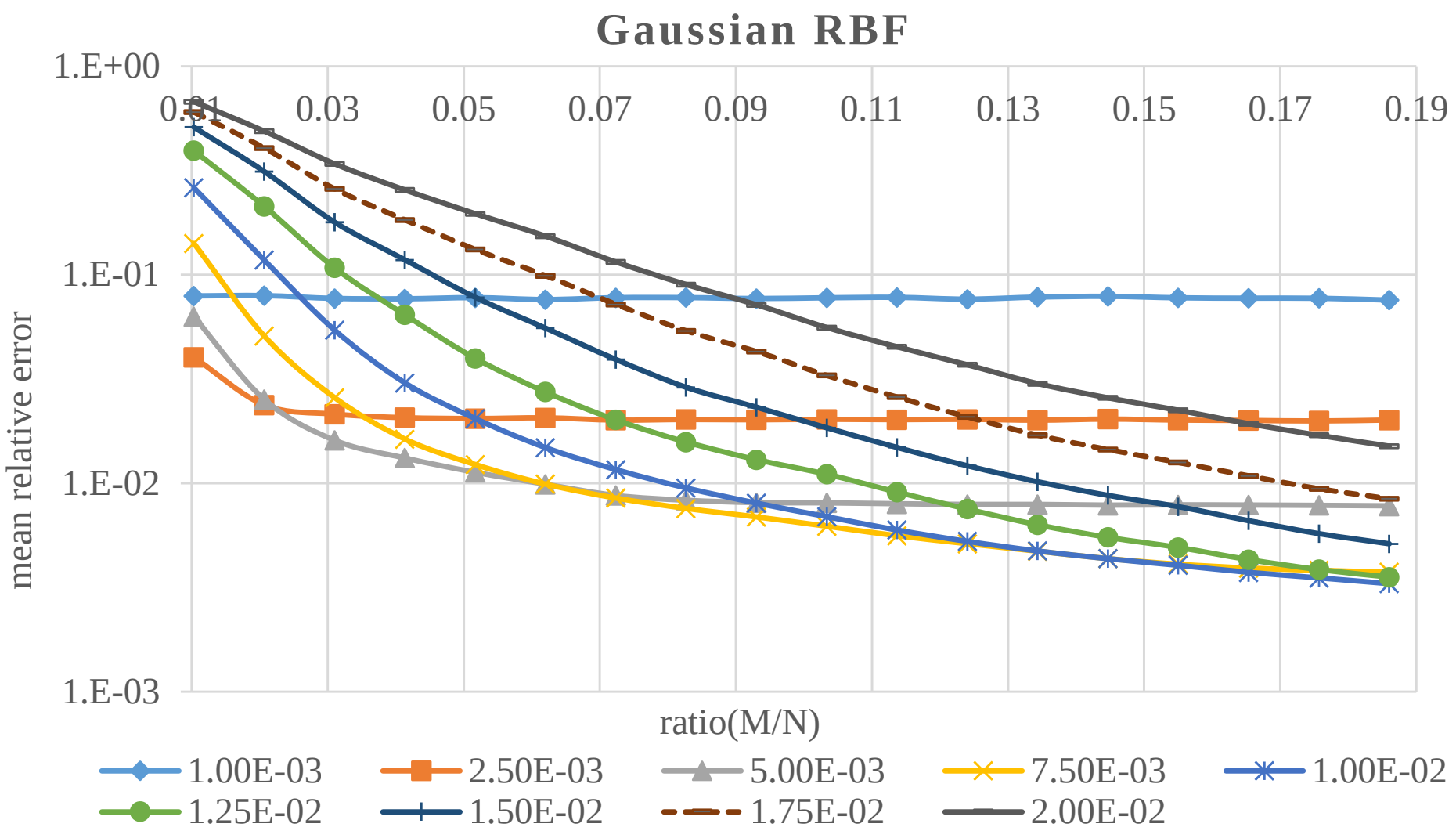

Figure 21: The mean relative error of the proposed RBF approximation of mount Veľký Rozsutec dataset for different numbers of reference points, Gaussian RBF with different shape parameters $\alpha$.

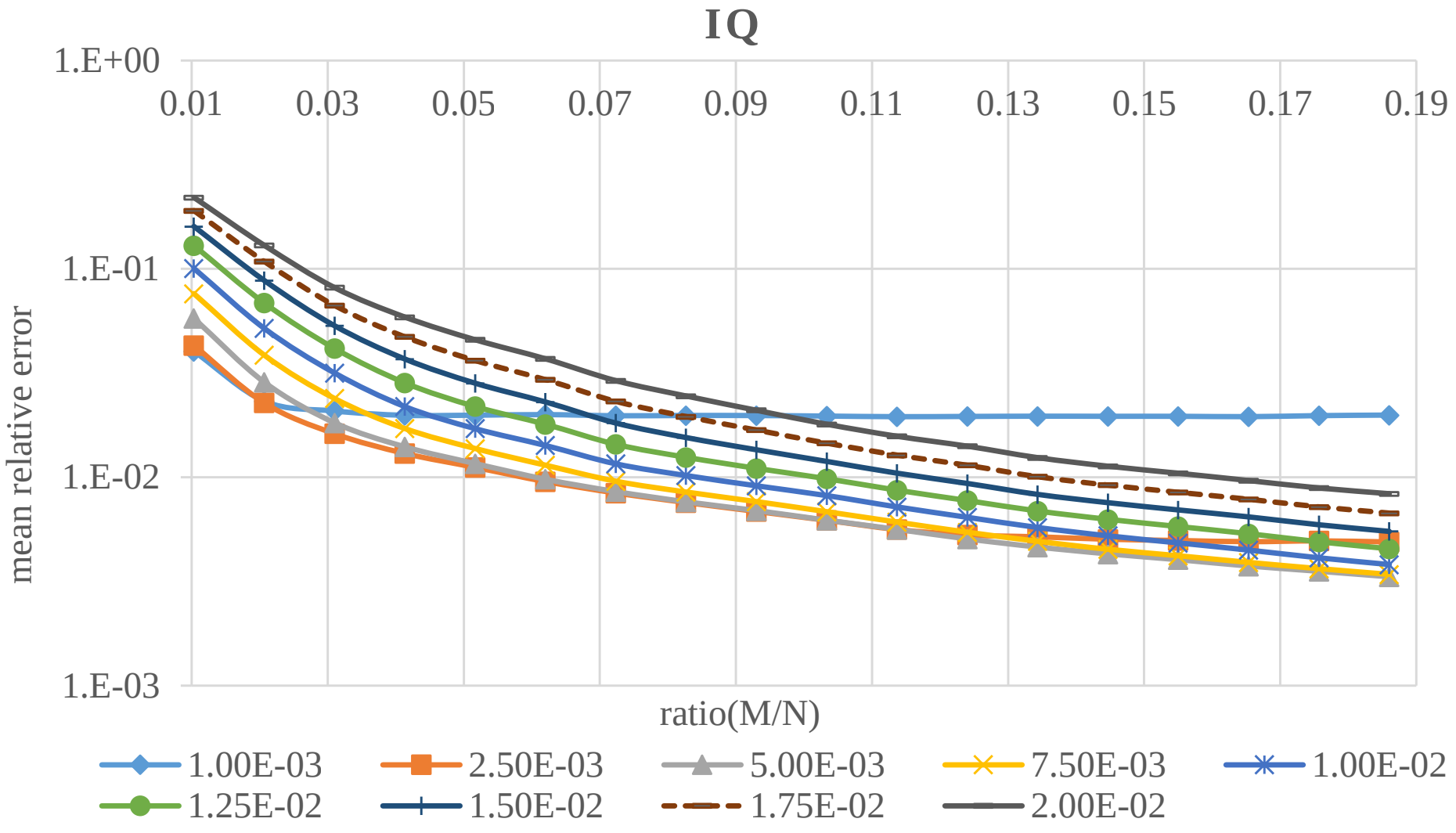

Figure 22: The mean relative error of the proposed RBF approximation of mount Veľký Rozsutec dataset for different numbers of reference points, IQ RBF with different shape parameters $\alpha$.

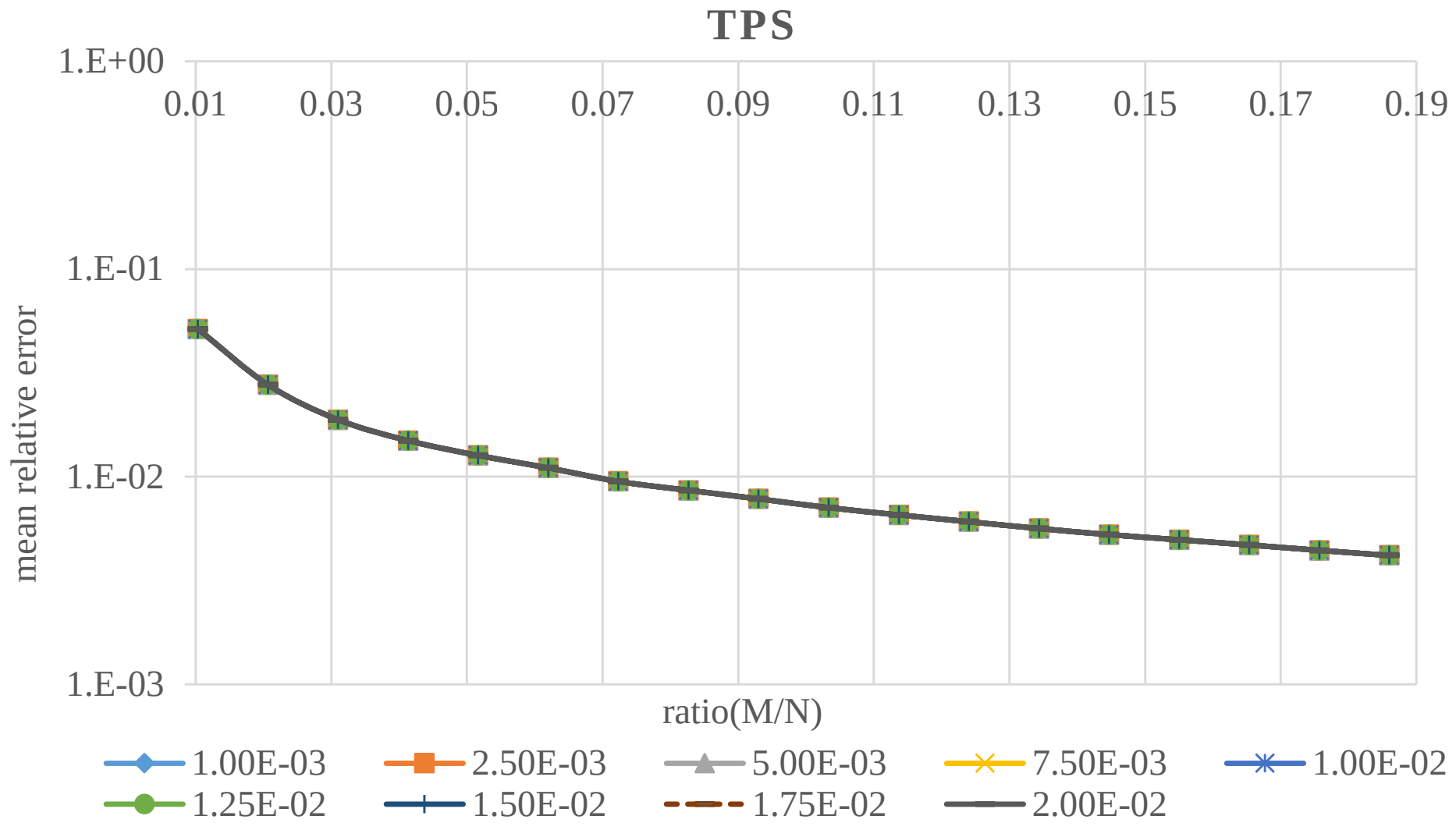

Figure 23: The mean relative error of the proposed RBF approximation of mount Veľký Rozsutec dataset for different numbers of reference points, TPS RBF with different shape parameters $\alpha$.

Finally, note that the results for RBF approximation with linear reproduction are very similar to the proposed RBF approximation.

## 6. Conclusion

Comparisons of different methods of RBF approximation with respect to various criteria were presented. The proposed RBF approximation introduced in Sect. 3 gives the best results due to the smallest error. The RBF approximation with a linear reproduction can be influenced by placement of the given dataset in space. Therefore, it is appropriate that the translation of the estimated center of gravity to the origin of the coordinate system is made as the first step. The worst results according to error were obtained using the RBF approximation using Lagrange multipliers. Moreover, this method of approximation has unpredictable behavior, the matrix for RBF approximation using Lagrange multipliers is mostly ill-conditioned and its size is high, i.e. it is of the $(M + N) \times (M + N)$ size.

The experiments proved that the proposed RBF approximation gives significantly better result over other methods used in the experiments described above. It also offers a possible data compression as the matrix is only $M \times M$, where $M \ll N$, which is a significant factor for large datasets processing. On the other hand, experiments made also proved that all methods have problems with the preservation of sharp edges if global functions are used.

Future work will be devoted to evaluation of Compactly-Supported RBFs (CS-RBFs) which will lead to sparse matrices, decrease of memory requirements and significant increase of speed of computation. A special attention will be given to finding optimal shape parameters which is critical for the RBF approximation quality.

## Acknowledgments

The authors would like to thank their colleagues at the University of West Bohemia, Plzen, for their discussions and suggestions, and anonymous reviewers for their valuable comments and hints provided. The research was supported by the National Science Foundation GAČR project GA17-05534S and partially supported by SGS-2016-013.